\documentstyle[amscd,amssymb,verbatim,12pt]{amsart}

\theoremstyle{plain}
\newtheorem{Thm}{Theorem}
\newtheorem{Cor}{Corollary}
\newtheorem{Lem}{Lemma}
\newtheorem{Prop}{Proposition}

\theoremstyle{definition}
\newtheorem{Def}{Definition}

\theoremstyle{remark}

\newtheorem{Rem}{Remark}

\numberwithin{equation}{section}
\renewcommand{\rm}{\normalshape}

\newif\ifShowLabels
\ShowLabelstrue
\newdimen\theight
\def\TeXref#1{%
	\leavevmode\vadjust{\setbox0=\hbox{{\tt
		\quad\quad  {\small \rm #1}}}%
	\theight=\ht0
	\advance\theight by \lineskip
	\kern -\theight \vbox to
	\theight{\rightline{\rlap{\box0}}%
	\vss}%
	}}%

\ShowLabelsfalse

\renewcommand{\sec}[2]{\section{#2}\label{S:#1}%
	\ifShowLabels \TeXref{{S:#1}} \fi}
\newcommand{\ssec}[2]{\subsection{#2}\label{SS:#1}%
	\ifShowLabels \TeXref{{SS:#1}} \fi}
\newcommand{\sssec}[2]{\subsubsection{#2}\label{SSS:#1}%
        \ifShowLabels \TeXref{{SSS:#1}} \fi}

\newif\ifShowLabels
\ShowLabelstrue
\newdimen\theight
\def\TeXref#1{%
	\leavevmode\vadjust{\setbox0=\hbox{{\tt
	\quad\quad  {\small \rm #1}}}%
	\theight=\ht0
	\advance\theight by \lineskip
	\kern -\theight \vbox to
	\theight{\rightline{\rlap{\box0}}%
	\vss}%
	}}%

\ShowLabelsfalse

\newcommand{\refs}[1]{\ref{S:#1}}
\newcommand{\refss}[1]{\ref{SS:#1}}
\newcommand{\refsss}[1]{\ref{SSS:#1}}
\newcommand{\refe}[1]{\eqref{E:#1}}
\newcommand{\reft}[1]{\ref{T:T:#1}}
\newcommand{\refl}[1]{\ref{L:#1}}
\newcommand{\refp}[1]{\ref{P:#1}}
\newcommand{\refc}[1]{\ref{C:#1}}
\newcommand{\refd}[1]{\ref{D:#1}}
\newcommand{\refr}[1]{\ref{R:#1}}

\newenvironment{thm}[1]%
	{ \begin{Thm} \label{T:#1}  \ifShowLabels \ref{T:T:#1} \fi }%
	{ \end{Thm} }

\newcommand{\th}[1]{\begin{thm}{T:#1} \sl}
\renewcommand{\eth}{\end{thm} }

\newenvironment{lemma}[1]%
	{ \begin{Lem} \label{L:#1}  \ifShowLabels \TeXref{L:#1} \fi }%
	{ \end{Lem} }
\newcommand{\lem}[1]{\begin{lemma}{#1} \sl}
\newcommand{\elem}{\end{lemma}}

\newenvironment{propos}[1]%
	{ \begin{Prop} \label{P:#1}  \ifShowLabels \TeXref{P:#1} \fi }%
	{ \end{Prop} }
\newcommand{\prop}[1]{\begin{propos}{#1} \sl}
\newcommand{\eprop}{\end{propos}}

\newenvironment{corol}[1]%
	{ \begin{Cor} \label{C:#1}  \ifShowLabels \TeXref{C:#1} \fi }%
	{ \end{Cor} }
\newcommand{\cor}[1]{\begin{corol}{#1}\sl }
\newcommand{\ecor}{\end{corol}}

\newenvironment{defeni}[1]%
	{ \begin{Def} \label{D:#1}  \ifShowLabels \TeXref{D:#1} \fi \rm}%
	{ \end{Def} }
\newcommand{\defe}[1]{\begin{defeni}{#1}\sl }
\newcommand{\edefe}{\end{defeni}}

\newenvironment{remark}[1]%
	{ \begin{Rem} \label{R:#1}  \ifShowLabels \TeXref{R:#1} \fi \rm }%
	{ \end{Rem} }
\newcommand{\rem}[1]{\begin{remark}{#1}}
\newcommand{\erem}{\end{remark}}

\newtheorem{remar}{Remark}

\newcommand{\rema}{\begin{remar}\rm}
\newcommand{\erema}{\end{remar}}

\newcommand{\eq}[1]%
	{ \ifShowLabels \TeXref{E:#1} \fi 
	   \begin{equation} \label{E:#1} }
\newcommand{\eeq}{ \end{equation} }

\newenvironment{proof}{\smallskip\noindent{\it Proof.\quad}}{\qquad$\square$
\par\smallskip}
\newcommand{\prf}{ \begin{proof} }
\newcommand{\epr}{ \end{proof}  }

\newcommand\alp{\alpha}		
\newcommand\bet{\beta}		
\newcommand\gam{\gamma}		\newcommand\Gam{\Gamma}
		\newcommand\Del{\Delta}
\newcommand\eps{\varepsilon}

\newcommand\kap{\kappa}
\newcommand\lam{\lambda}		\newcommand\Lam{\Lambda}

\newcommand\ome{\omega}		

\newcommand\calA{{\cal{A}}}
\newcommand\calB{{\cal{B}}}
\newcommand\calC{{\cal{C}}}
\newcommand\calD{{\cal{D}}}

\newcommand\calF{{\cal{F}}}
\newcommand\calG{{\cal{G}}}

\newcommand\calK{{\cal{K}}}
\newcommand\calL{{\cal{L}}}

\newcommand\calO{{\cal{O}}}

\newcommand\calR{{\cal{R}}}
\newcommand\calS{{\cal{S}}}

\newcommand\QQ{\Bbb{Q}}

\newcommand\PP{\Bbb{P}}
\renewcommand\AA{\Bbb{A}}

\newcommand\FF{\Bbb{F}}
\newcommand\GG{\Bbb{G}}

\newcommand\ZZ{\Bbb{Z}}

\newcommand\NN{\Bbb{N}}

\newcommand\sdp{\times \hskip -0.3em {\raise 0.3ex
\hbox{$\scriptscriptstyle |$}}} 

\newcommand\cf{{\rm \,cf\,}}

\newcommand\End{\operatorname{End\,}}
\newcommand\Ext{\operatorname{Ext}}

\newcommand\rank{\operatorname{rank}}

\newcommand\Hom{\operatorname{Hom}}
\newcommand\RHom{\operatorname{RHom}}

\newcommand\Ker{\operatorname{Ker}}
\newcommand\Perv{\operatorname{Perv}}

\newcommand\pr{\text{pr}}

\newcommand\tr{\operatorname{tr}}
\newcommand\Tr{\operatorname{Tr}}

\newcommand\oA{{\overline{A}}}

\newcommand\oX{{\overline{X}}}

\newcommand\oalp{{\overline{\alpha}}}

\newcommand\okap{{\overline{\kap}}}

\newcommand\opi{{\overline{\pi}}}

\newcommand\tilA{{\widetilde{A}}}

\newcommand\tilE{{\widetilde{E}}}

\newcommand\tilG{{\widetilde{G}}}

\newcommand\tilj{{\widetilde{j}}}
\newcommand\tilm{{\widetilde{m}}}

\newcommand\tilp{{\widetilde{p}}}

\newcommand\tilY{{\widetilde{Y}}}

\newcommand\tilDel{{\widetilde{\Delta}}}

\newcommand\tilphi{{\widetilde{\phi}}}
\newcommand\tilpi{{\widetilde{\pi}}}

\newcommand\tilPhi{{\widetilde{\Phi}}}

\newcommand\x{\times }
\newcommand\ten{\otimes}

\renewcommand{\>}{\rangle}
\newcommand{\<}{\langle}
\newcommand\la{\langle}
\newcommand\ra{\rangle}

\newcommand{\lan}{\langle}
\newcommand{\ran}{\rangle}

\newcommand{\Id}{\operatorname{Id}}

\newcommand{\und}{\underline}

\renewcommand{\Im}{\operatorname{Im}}
\newcommand{\Ga}{\Gamma}
\newcommand{\ga}{\gamma}

\newcommand{\ed}{\qed\vspace{3mm}}

\newcommand{\Spec}{\operatorname{Spec}}

\newcommand{\ov}{\overline}

\newenvironment{exs}{\vspace{3mm} \noindent {\bf
Examples.}}{\vspace{3mm}}
\theoremstyle{definition}

\newcommand{\Pf}{\noindent {\it Proof}}

\renewcommand{\a}{\alpha}
\renewcommand{\b}{\beta}

\newcommand{\wt}{\widetilde}

\newcommand{\A}{{\Bbb A}}

\newcommand{\Q}{{\Bbb Q}}

\newcommand{\id}{\operatorname{id}}

\newcommand{\D}{{\cal D}}
\newcommand{\De}{\Delta}

\newcommand{\lrar}[1]{\begin{picture}(50,10)(-25,-5)
\put(-25,0){\vector(1,0){50}}
\put(0,5){\makebox(0,0)[b]{\mbox{$#1$}}}
\end{picture}}

\newcommand{\ldar}[1]{\begin{picture}(10,50)(-5,-25)
\put(0,25){\vector(0,-1){50}}
\put(5,0){\mbox{$#1$}}
\end{picture}}

\newcommand{\qlbg}{{\overline \QQ}_{\ell,G}}

\newcommand\dspace{\lineskip=2pt\baselineskip=18pt\lineskiplimit=0pt}
\dspace
\newcommand\fq{\FF_q}
\newcommand\fqb{{\overline \FF_q}}
\newcommand\gfq{G(\fq)}
\newcommand\fr{\text{Fr}}
\newcommand\qlb{{\overline \QQ_{\ell}}}
\newcommand\chl{\calK_{\calL}}
\newcommand\al{\calA_{\calL}}

\newcommand\Fo{\text{\bf Four}} 

\newcommand\Ra{\text{\bf Rad}}

\newcommand\frw{\text{Fr}_w}
\newcommand\avg{\text{{\bf Av}}_G}
\newcommand\wl{W_{\calL}}
\newcommand\wll{W_{\calL_1,\calL_2}}
\newcommand\wtheta{W_{\theta}}
\newcommand\mon{\text{mon}}
\newcommand\reg{\text{reg}}
\newcommand\brw{\text{Br}_W}

\begin{document}

\dspace
\title[Kazhdan-Laumon representations]
{Kazhdan-Laumon representations of Chevalley groups,
character sheaves and some generalization of the Lefschetz-Verdier 
trace formula}

\author[Alexander Braverman and Alexander Polishchuk]{\  }

\address{A.~B.: room 2-175, MIT, Department of Mathematics, 77 
Mass. Ave., Cambridge, MA, 02139, USA
\newline
A.~P. Department of Mathematics, Harvard University, 1 Oxford St., Cambridge,
MA, 02139, USA}
\maketitle
\centerline{Alexander Braverman
and Alexander Polishchuk
\footnote{The work of both authors was partially supported by the
Natural Scince Foundation}}

\begin{abstract}
Let $G$ be a reductive group over a finite field $\fq$ and let $\gfq$
be its group of $\fq$-rational points. In 1976 P.~Deligne and G.~Lusztig
(following a suggestion of V.~Drinfeld for $G=SL(2)$) associated
to any maximal torus $T$ in $G$ defined over $\fq$
certain remarkable series (virtual) representations of the group $G(\fq)$.
These representations were realized in \'e tale $\ell$-adic cohomology
of certain algebraic varieties over $\fq$, (called the Deligne-Lusztig  
varieties). Using these representations G.~Lusztig has later given
a complete classification of representations of $\gfq$. He has also discovered
that characters of these representations might be computed from certain
geometrically defined $\ell$-adic perverse sheaves on the group $G$, called
character sheaves.

Despite the fact that both Deligne-Lusztig representations and character
sheaves are defined by some simple (and beautiful) 
algebro-geometric procedures,
the proof of the relation between the two is rather complicated (apart
from the case, when $T$ is split, where it becomes just the usual formula
for the character of an induced representation).
In 1987 D.~Kazhdan and G.~Laumon proposed a different (conjectural)
geometric way 
for constructing representations of $\gfq$. Their idea was based
on exploiting the generalized Fourier-Deligne on the basic affine
space $X=G/U$ of $G$ (here $U$ is a maximal unipotent subgroup
of $G$).
In the first part of this paper we give a
rigorous definition of (almost all) Kazhdan-Laumon representations.
Namely for every $\fq$-rational maximal torus $T$ in $G$ and every
quasi-regular ({\it non-singular} in the terminology of \cite{dl})
character $\theta:\ T(\fq)\to\qlb^*$ we construct certain representation
$V_{\theta}$ of $\gfq$. We show that this representation is finite-dimensional
and that it is irreducible if $\theta$ is regular ({\it in a generic position}
in the terminology of \cite{dl}). Moreover, it follows from the proof of
irreducibility in essentially tautological way, that the character of
$V_{\theta}$ is given by the trace function of the corresponding Lusztig's
character sheaf. Roughly speaking, all the above is achieved by 
replacing the generalized Fourier-Deligne transformations on $X$ by
the suitable version of Radon transformations.

It follows from the above computation of the character of $V_{\theta}$, 
that $V_{\theta}$ is equivalent to the corresponding Deligne-Lusztig
representation $R_{\theta}$. The second part of this paper is devoted
to the construction of an explicit isomorphism between the two.
This is done by using certain generalization of the Lefschetz-Verdier
trace formula for Radon transformations, which we cannot establish
in the full generality (however, we believe that this is a technical
difficulty, which is there only because of authors' laziness).

\end{abstract}

\tableofcontents

\sec{}{Introduction}

\ssec{}{Preliminaries}Let $p$ be an prime number, $q$ --
a power of $p$ and let $\FF_q$ denote the finite field with $q$
elements. For any scheme $X$ defined over $\FF_q$ we denote by
$X(\fq)$ the set of all $\fq$ -rational points of of $X$. Let
also $\fr :\ X\to X$ denote the geometric Frobenius morphism.
We will also fix a prime number $\ell\neq p$. 

We will denote by $\calD(X)$ the bounded derived category of $\ell$-adic
constructible sheaves on $X$ and by $\Perv(X)$ the corresponding
category of perverse sheaves.
\ssec{dl}{Deligne-Lusztig representations}
Let $G$ be an algebraic reductive connected
group defined over $\fq$ (in fact, for the sake of simplicity, we will speak
mostly about semisimple simply connected groups; however, 
everything done in this paper goes through for any reductive group over
$\fq$ with only minor change)s. The representations of the finite group
$\gfq$
were classified by G.~Lusztig. As his main tool, G.~Lusztig used
an earlier construction of some amount of (virtual) representations of 
$\gfq$ by P.~Deligne and G.~Lusztig. Let us recall some facts about this
construction. 

Recall (cf. \cite{dl}) that one can associate to $G$ its abstract Cartan 
group $T$, together with canonical $\fq$-rational structure on it and a
root system (with the chosen system of positive roots). Let $W$ be the Weyl
group of $G$ which will be considered as a subgroup of the group
Aut$(T)$ of automorphisms of $T$. Any $w\in W^{\fr}$ defines a new $\fq$-
rational structure on $T$ by composing the old action of $\fr$ with
$w$. We will denote this  Frobenius morphism by $\frw$. 
Let $T_w(\fq)$ denote the group of $\fq$-rational points of $T$ 
with the above $\fq$-structure. Thus $T_w(\fq)=T(\fqb)^{\frw}$.

For any $w\in W$ Deligne and Lusztig have defined certain
algebraic variety $X_w$ over $\fqb$ endowed with an action
of the group $\gfq\times T_w(\fq)$. Thus we may
consider the virtual representation $\sum (-1)^iH^i_c(X_w,\qlb)$
of $\gfq$ and decompose it with respect to characters of $T_w(\fq)$.
For any character $\theta:\ T_w(\fq)\to \qlb^*$ we will denote
by $R_{{\theta},w}$ the corresponding virtual representation
of $G(\fq)$. In their paper Deligne and Lusztig showed
the following
\th{} Suppose that $\theta$ is regular (cf. \refss{characters}; note, that 
what we call regular character here is called {\it a character in a generic
position} in \cite{dl}). Then $R_{\theta, w}$
is an irreducible representation of $\gfq$ (in particular, it is a genuine
representation, not a virtual one).
\eth
\ssec{char-sheaves}{Character sheaves}Let $\calL$ be a one-dimensional tame
(\cf. \refss{characters}) 
$\ell$-adic local system on $T$. Suppose that we have
also fixed an isomorphism $\calL\simeq \frw^*\calL$ for some $w\in W$. 
In \cite{lu-char} G.~Lusztig
has defined certain Ad$G$-equivariant perverse sheaf $\chl$ on the
group $G$, equipped with the structure of a Weil sheaf over
$\fq$ (in the notations of \cite{lu-char} this is the sheaf $\chl^e$, where
$e$ denotes the unit element in $W$). Let $\tr(\chl)$ denote the 
``trace function'' of $\chl$
(by the Grothendieck faicseaux-fonctions correspondence).
This is an Ad$\gfq$-equivariant function on $\gfq$ with values
in $\qlb$. The following result is due to G.~Lusztig (cf., for example,
\cite{lu-90}).

\th{}
One has 
\eq{char-dl}
\tr(\chl)=\chi(R_{\theta, w})
\end{equation}
where $\chi(R_{\theta, w)}$ is the character of the virtual representation
$R_{\theta,w}$.
\eth
\ssec{kl}{Kazhdan-Laumon representations}Despite the fact that both sides
of the equality \refe{char-dl} are defined by simple geometric procedures,
the proof of \refe{char-dl} is not as straightforward as one would 
expect it (except for the
case $w=1$ where it becomes just the usual formula for the character of an
induced representation).
This phenomenon might have 2 explanations.

1) The Deligne-Lusztig varieties $X_w$ change as we extend the base field
$\FF_q$, which makes difficult the analysis of the behaviour of $R_{\theta,w}$
under extensions of the base field.

2) The Weil structure on the sheaf $\chl$ is defined in a rather inexplicit
way (it
is defined explicitly on the set $G_{rs}$ of regular semisimple
elements in $G$, and then one defines it on the whole
of $G$ by making use of the fact that $\chl$ is an intermediate extension
of its restriction to $G_{rs}$).

Let now $\calL$ be a tame quasi-regular one-dimensional local system on $T$
(this means that $\alp^*\calL$ is not isomorphic to the constant sheaf for
any coroot $\alp:\GG_m\to T$) endowed with an
isomorphism $\calL\simeq \frw^*\calL$.
One of the purposes of this paper is to give a simple geometric proof
of the fact that $\tr(\chl)$ is a character of some irreducible 
representation of $\gfq$. To do that we use a different construction
of representations of $\gfq$ suggested by D.~Kazhdan and G.~Laumon in
\cite{kl}. We also give a different way to define the structure of a
Weil sheaf on $\chl$. Both these steps are done by exploiting the
generalized Fourier-Deligne transformations on the basic affine space
$X=G/U$ ($U$ is a maximal unipotent subgroup of $G$).

Let us give a brief recollection of the variant of the Kazhdan-Laumon
construction which will be used in this paper. We would like to imitate
``as much as possible'' the construction of principle series. Recall
that the latter are defined as follows.
Let $X$ be again the basic affine space of $G$. 
This is a quasi-affine algebraic variety defined
over $\fq$. It admits a natural action of the group $G\times T$ ($T$ acts
there since we may choose a maximal torus in $G$ which normalizes $U$). 
Therefore,
the set $X(\fq)$ admits a natural action of $\gfq\x T(\fq)$. Hence, given
a character $\theta$ of $T(\fq)$ we may construct a representation $V(\theta)$
of $G(\fq)$ (this is just the representation of $\gfq$ in the space of
functions on $X(\fq)$ which change according to character $\theta$ under
the action of $T(\fq)$).

Now, we would like to construct series of representations 
which are parametrised
by characters of other maximal tori in $G$ which are defined over $\fq$. 
According to \cite{dl}
different conjugacy classes of maximal tori in $G$, which are defined over
$\fq$ are parametrised by twisted conjugacy classes in $W$ (if $G$ is split, 
then twisted conjugacy classes are just the same as the usual conjugacy 
classes; in the general case cf. \refss{characters} for the definition).

In \cite{kl} 
Kazhdan and Laumon have described certain idea how this can be done. 
Unfortunately, their construction relied on some conjectures from 
homological algebra which are still unknown. However, one can use
certain modification of the Kazhdan-Laumon construction, which is already
well-defined.

It is observed in \cite{kl} that the category $\Perv(X)$ of $\ell$-adic
perverse sheaves on $X$ admits a natural action of the braid group, 
corresponding
to $W$ (cf. \refss{four-on-x}). 
In \refss{four-on-x} we define certain subcategory
$\Perv^0(X)$ of $\Perv(X)$, on which the above action of the braid 
group can be reduced to an action of $W$. This subcategory is ``big enough''
(in particular, it is still invariant under the action of $G\x T$).

Consider the category $\Perv^0_{w}$ of ``$w$-Weil sheaves'' on $X$.
An object of $\Perv^0_{w}$ is an object $A$ of $\Perv^0(X)$
together with an isomorphism $\alp: A\simeq \fr^*\Phi_w(A)$ where 
$\Phi_w$ is the functor on $\Perv^0(X)$ defining the action of $W$ on it
and $\fr :X\to X$ is the corresponding Frobenius morphism.
For $w=1$ we get just the usual notion of Weil sheaves. 
Consider now the space $K:=K(\Perv^0(X))\ten \qlb$ ($K(\Perv^0(X))$
is the Grothendieck group of $\Perv^0(X)$). This space is an
infinite-dimensional representation
of the group $\gfq\x T_w(\fq)$. Next, one can define certain quotient 
$V_{w}$ of this space. The quotient is defined in such a way that for $w=1$ if
we replace $\Perv^0(X)$  by Perv$(X)$ then we get just the
space of functions on the finite set $X(\fq)$. For a character $\theta$
of $T_w(\fq)$ we denote $V_{\theta,w}$ the corresponding representation
of $\gfq$.

\th{} Let $\theta$ be any  character of $T_w(\fq)$.
Then the Kazhdan-Laumon representation $V_{\theta,w}$ is  finite-dimensional
and nonzero. If $\theta$ is regular ( {\it in a generic position} in
the terminology of \cite{dl}) then $V_{\theta,w}$
is irreducible. Moreover, if $\theta$ is non-singular (in the 
sense of \cite{dl}) then the character of $V_{\theta,w}$ is equal
to $\tr(\chl)$, where $\calL$ is the one-dimensional local system on
$T$, which corresponds to $\theta$. 
The last fact is established by a simple geometric argument which is
essentially ``not different'' from the case of principal series).  
\eth

One can also write down an explicit geometric isomorphism between
$V_{\theta, w}$   and the corresponding Deligne-Lusztig representation
(this is done in Section \refs{delus}). This construction 
requires certain generalization  of the Lefschetz-Verdier 
trace formula, which is the subject of Section \refs{trace}.  
In this way we obtain a new proof of the fact that characters of
the Deligne-Lusztig representations can be computed using the corresponding
character sheaves (modulo the assumption that our generalized trace fromula
holds in the case we need).
\ssec{}{Contents}This text is organized as follows. In section
\refs{fourier} we establish some facts about the Fourier-Deligne
transforms on the basic affine space, which will be needed in the
sequel. In Section \refs{char} we recall some basic facts about the
sheaf $\chl$ and give a definition of the Weil structure on $\chl$
which is slightly different from that of G.~Lusztig and uses 
the Fourier-Deligne transforms on $X$. In Section 
\refs{kl} we give the precise 
definition of Kazhdan-Laumon representations and state our main results
about them. Section \refs{proof} is devote to the proof of these
results. 
Finally, in Section \refs{trace} we discuss certain generalization of
the Lefschetz-Verdier trace formula and in Section \refs{delus} use
this generalization in order to construct a geometric isomorphism
between Kazhdan-Laumon and Deligne-Lusztig representations.

\ssec{}{Acknowledgements}We would like to thank D.~Kazhdan for
posing the problem and a lot of very useful discussions. 
Some of the results of this paper were presented in PhD thesis of the
first author, who would like to thank his advisor J.~Bernstein for
continuous interest in this work and many helpul advices. 
We are also grateful to A.~Beilinson, R.~Bezrukavnikov, D.~Gaitsgory, 
G.~Laumon and G.~Lusztig for useful discussions and their remarks about
early versions of this text.

The work of both authors was partially supported by the
National Science Foundation. 
\sec{fourier}{Generalized Fourier-Deligne and Radon transformations
on the basic affine space}

\ssec{}{Fourier-Deligne and Radon transformations on a vector bundle}
\sssec{four-bundle}{Definition of Fourier transform}
Let $S$ be a scheme of finite type over $k=\fqb$ 
and let $\pi_E:E\to S$ be a vector
bundle over $S$ of rank $r$ and let 
$\pi_{E^{\vee}}:E^{\vee}\to S$ be the dual bundle. 
Fix a nontrivial character $\psi$ of $\fq$ and denote by
$\Fo_{\psi}:\calD(E)\to \calD(E^{\vee})$ the Fourier-Deligne transform (cf. 
\cite{KatzLaum}). Recall that it is defined as follows. Let $\calL_{\psi}$ 
denote the Artin-Schreier sheaf on $\AA^1=\AA^1_k$. Let also 
$\mu:~E\x_S E^{\vee}\to \AA^1$ denote the pairing map and 
let $\pr_1,\pr_2$ denote the
projections from $E\x_S E^{\vee}$ to the first and to the second variable
respectively. Then, by definition
\eq{}
\Fo_{\psi}(A)=\pr_{2!}(\pr_1^*A\ten \mu^*\calL_{\psi})[r]
\end{equation}
for any $A\in \calD(E)$. We will omit the subscript ``$\psi$'' when it 
does not lead to a confusion.

Let $inv:~\GG_m\to \GG_m$ denote the inverseion map 
(i.e. $inv(\lam)=\lam^{-1}$
and let $m_E:~\GG_m\x E\to E$ denote the multiplication map (defined by 
$m_E(\lam,e)=\lam e$). Then it is easy to see that one has a canonical
isomorphism of functors
\eq{}
(\id\x \Fo)\circ m_E^*\simeq  (inv\x \id)^*\circ m_{E^{\vee}}^*\circ \Fo
\end{equation}
 
\sssec{rad-bundle}{Definition of Radon transform}
Let us now define the version of the Radon transformation that will be needed
later. The notations will be as above. For a vector bundle $\pi_E:E\to X$
let us denote by $j_E:\tilE\to E$ the embedding of the complement to
the zero section on $E$. Let now $Z_E\subset E\times E^{\vee}$ be the
closed subvariety in $E\times E^{\vee}$ defined by
\eq{}
Z_E=\{ (e,e^{\vee})\in E\times E^{\vee}|\ \la e,e^{\vee}\ra = 1\}
\end{equation}
Note that $Z$ lies, in fact, inside $\tilE\x \tilE^{\vee}$. 
Let $p_1,p_2$ denote respectively the  projections from
$Z$ to $\tilE$ and $\tilE^{\vee}$.
We now define the Radon transform $\Ra:\ \calD(\tilE)\to 
\calD(\tilE^{\vee})$ by
\eq{}
\Ra(A)=p_{2!}p_1^*(E)[r-1]
\end{equation}
\sssec{dzero}{The categories $\calD^{mon}(E)$, $\calD^0(E)$ and 
$\calD^{reg}(E)$}
The main difference between the Fourier-Deligne transform 
and 
the Radon transform, is that the latter 
does not map in general perverse sheaves into perverse
sheaves. However, below we will define certain full subcategory of 
$\calD(E)$ (resp. $\calD(E^{\vee})$), which we will denote by
$\calD^{0}(E)$ (resp. $\calD^{0}(E^{\vee})$), such that the 
restrictions of
the functors $\Fo$ and $\Ra$ to $\calD^{0}(E)$ will be canonically
isomorphic. This, in particular, will imply that $\Ra$ maps perverse
objects, lying in $\calD^{0}(E)$, into perverse ones. 

First of all, let us define the category of monodromic sheaves.
\defe{mon}
1) A complex $A\in \calD(E)$ is called {\it monodromic} if there exists
$n\in \NN$ such that for every $i\in \ZZ$, $H^i(A)$ has a filtration, whose
factors are equivariant with respect to the action of
the group $\GG_m$ on $E$, defined by $\lam(e)=\lam^ne$ for 
$\lam\in\GG_m$ and $e\in E$ (here $H^i(A)$ denotes the $i$-th perverse
cohomology of $A$). We denote by $\calD^{\mon}(E)$ the
full subcategory of $\calD(E)$ consisting of monodromic complexes.

2) A monodromic complex $A$ is called regular if for every $i\in \ZZ$,
$H^i(A)$ has a filtration, whose quotients are $(\GG_m,\calL)$-equivariant,
where $\calL$ is some non-constant local system on $\GG_m$. We denote by
$\calD^{\reg}(E)$ the full subcategory of $\calD(E)$, consisting of
regular monodromic complexes. 
\edefe

\rem{}
The word ``equivariant'' in Definition \refd{mon} means, in fact,
``can be given an equivariant structure''.
\erem

\rem{monod}Part 1 of Definition \refd{mon} makes sense for sheaves on arbitrary
variety $X$, endowed with an action of an algebraic torus $T$.
\erem

\rem{}Let $X$ be any variety and let $A\in \calD(X)$. Supose that
for every $i\in \ZZ$ we are given a filtration of $H^i(A)$ with
certain quotients. Then we will say that $A$ {\it is glued}
from those quotients.
\erem

We will define now a third category $\calD^0(E)$, which will
be intermediate between $\calD^{\reg}(E)$ and $\calD^{\mon}(E)$ (i.e.
we will have natural inclusions $\calD^{\reg}(E)\subset \calD^0(E)
\subset \calD^{\mon}(E)$ of full subcategories).

Let $i_E:\ X\to E$ denote the embedding of the zero section.
By definition, the category $\calD^0(E)$ is a full subcategory of 
$\calD^{\mon}(E)$,
consisting of objects $A$ of $\calD^{\mon}(E)$, such that

1) $\pi_{E*}A=\pi_{E!}A=0$

2) $i_E^*A=i_E^!A=0$

Note, that it follows from the definition that $\calD^0(E)$ can be
identified with the full subcategory of $\calD^{\mon}(\tilE)$, consisting
of those $A\in \calD^{\mon}(\tilE)$ for which the canonical morphism
$j_{E!}A\to j_{E*}A$ is an isomorphism and  
$\pi_{E*}(j_{E*}A)=\pi_{E!}(j_{E!}A)=0$. Therefore, we will sometimes write
$\calD^0(\tilE)$ instead of $\calD^0(E)$, in order to 
emphasize, that the objects of this category may be regarded as sheaves on
$\tilE$.

We will denote by $\Perv^{\mon}(E)$, $\Perv^0(E)$ and $\Perv^{\reg}(E)$
the corresponding categories of perverse sheaves.
\th{rad-four}   
 
1) Both functors $\Ra$ and $\Fo$ map $\calD^{\mon}(E)$ to
$\calD^{\mon}(E^{\vee})$, $\calD^{\reg}(E)$ to
$\calD^{\reg}(E^{\vee})$ and $\calD^0(E)$ to
$\calD^0(E^{\vee})$ (more precisely, one should say that
$\Ra$ maps $\calD^{\mon}(\tilE)$ to $\calD^{\mon}(\tilE^{\vee})$).

2) One has canonical isomorphism 
\eq{}
\Ra\simeq j^*_E\circ \Fo \circ j_{E!}
\end{equation}
of functors, going from $\calD^{\mon}(\tilE)$ to 
$\calD^{\mon}(\tilE^{\vee})$.
\eth

\prf The proof of point (1) of Theorem \reft{rad-four} is completely
straightforward and it is left to the reader. The proof of
(2) is essentially a repetition of the proof of Theorem 9.13 in
\cite{br}. Let us, however, present it for the sake of completeness.

The proof follows from the following result.
\lem{technical} 
Let $Y$ be a scheme over $\fqb$. Let $K\in \calD^{\mon}(\GG_m\x Y)$ and 
let $\tau:\GG_m\hookrightarrow \AA^1$
denote the natural embedding. Then one has canonical isomorphism 
\eq{}
\pr_!((\tau\x \id)_!K\ten \calL_{\psi})\simeq i_1^*K[-1]
\end{equation}
where $i_1$ denotes the embedding of $1\x Y$ in $\GG_m\x Y$.
\elem

\noindent
The proof is left to the reader.

Let us now show how the lemma implies point (2) of the theorem.
Consider the following diagram:
\begin{center}
  \begin{picture}(9,7)

	\put(3,0){\makebox(1,1){$E$}}
	\put(0,3){\makebox(1,1){$\tilE$}}
	\put(3,6){\makebox(1,1){$\tilE\underset{S}\times E^{\vee}$}}
	\put(6,3){\makebox(1,1){$E\underset{S}\times E^{\vee}$}}
	\put(9,0){\makebox(1,1){$E^{\vee}$}}

	\put(1,3){\vector(1,-1){2}}
	\put(4,6){\vector(1,-1){2}}
	\put(7,3){\vector(1,-1){2}}
	\put(3,6){\vector(-1,-1){2}}
	\put(6,3){\vector(-1,-1){2}}

	\put(1,1){\makebox(1,1){$j$}}
	\put(1,5){\makebox(1,1){$\pr_1$}}
	\put(5,5){\makebox(1,1){$j\x \id$}}
	\put(8,2){\makebox(1,1){$\pr_2$}}
	\put(5,1){\makebox(1,1){$\pr_1$}}

  \end{picture}
\end{center}

Let ${\widetilde \pr}_2$ denote the composition of $\pr_2$ with $j\x \id$.
Then it follows that
\eq{}
\Fo(j_!A)\simeq  {\widetilde \pr}_{2!}(\pr_1^*A\ten 
(j\x \id)^*\mu^*\calL_{\psi})[r]
\end{equation}

On the other hand, ${\widetilde \pr}_{2}=\pr\circ (\mu,\pr_2)$, where 
$\mu\ :\tilE\x _S E^{\vee}\to \AA^1$ is the pairing map and
$\pr:\ \AA^1\x E^{\vee}$ is the projection to the second multiple.
Then
\begin{align}
\Fo&(j_!A)\simeq  (\pr\circ (\mu,\pr_2))_!(\pr_1^*A\ten 
(j\x \id)^*\mu^*\calL_{\psi})[r]\simeq \\
&\pr_!
((\mu,\pr_2)_!\pr_1^*A\ten (\calL_{\psi}\boxtimes 
{\overline \QQ}_{l,E^{\vee}}))[r]
\simeq
i_1^*(\mu,\pr_2)_!\pr_1^*A[r-1]
\end{align}
where the last two isomorphisms are given respectively by the projection
formula and Lemma \refl{technical}. On the other hand (by base change)
$i_1^*(\mu\x \id)_!\pr_1^*A[r-1]$ is isomorphic to 
$p_{2!}p^*_1A[r-1]$, which finishes the proof.
\epr

\cor{}Let $A\in \Perv^{\mon}(\tilE)$ and suppose that $j_{E!}A$ is
perverse. Then $\Ra(A)$ is perverse. In particular, $\Ra$ maps
$\Perv^0(E)$ to $\Perv^0(E)$. 
\ecor

\noindent
{\bf Warning 1.}\quad It follows from the definition, that both
functors $\Fo$ and $\Ra$ ``commute'' with the functor $\fr^*$ of 
inverse image with respect to the Frobenius morphism. This means 
that we have natural isomorphisms of functors 
$\nu_F:~\Fo\circ\fr^*{\widetilde \to} \fr^*\circ\Fo$ 
and $\nu_R:~\Ra\circ\fr^* {\widetilde \to} \circ\fr^*\Ra$. During the proof
of Theorem \reft{rad-four} we have constructed an isomorphism 
between the restrictions of $\Ra$ and $\Fo$ on the category
$\calD^0(E)$. However, this isomorphism does not transform 
$\nu_F$ into $\nu_R$ (even if we change $\nu_F$ or $\nu_R$ by a Tate twist). 
This is the reason why Radon transform (in the above
sense) is not equal to Fourier transform, when one passes from sheaves
to functions on $E(\fq)$ (in fact, the difference between Radon and Fourier
transform on the level of functions is given by certain $\Gam$-function,
attached to the field $\fq$).

\noindent
{\bf Warning 2.}\quad The term ``Radon transform'' is usually used
in the literature for a little different transformation, going
from $\calD(\PP(E))$ to $\calD(\PP(E^{\vee}))$, where $\PP(E)$ and
$\PP(E^{\vee})$ denote the corresponding projectivized bundles
(cf. \cite{br}, for example). We would like to note that if 
one identifies $\calD(\PP(E))$ with the derived 
category of $\GG_m$-equivariant sheaves on $\tilE$, then our Radon
transformation does not coincide with that of \cite{br}. For example,
point 2 of Theorem \reft{rad-four} for the Radon transform discussed
in \cite{br} does not hold in general, without some additional hypopaper
on $A$ (it does hold, however, for $\GG_m$-equivariant $A\in \calD^0(E)$). 
\ssec{symplectic}{Fourier transform on a symplectic vector bundle}Suppose now
that our vector bundle $\pi_E:~E\to S$ is endowed with some fiberwise
symplectic structure $\ome$. Then, using $\ome$ we can identify 
$E$ with $E^{\vee}$ and hence, we may consider Fourier transform as
a functor from $\calD(E)$ to itself. Let us denote this functor 
simply by $\calF(n)$. I.e. $\calF$ is the symplectic
Fourier-Deligne transform with a Tate twist by $-n$ 
(of course, one should remember that 
$\calF$ depends also on $\psi$). The following
result is due to A.~Polishchuk (\cf. \cite{po}). 
This result will not be used in the sequel, but it is useful to
keep it in mind.
\prop{}
There exists canonical isomorphism of functors $c:~\calF^2{\widetilde \to}
\text{Id}$. The pair $(\calF(n),c)$ defines an action of the group
$\ZZ_2$ on the category $\calD(E)$ (cf. \cite{De1} for this notion),
i.e. $c$ satisfies the following associativity condition:
\eq{}
\text{the two morphisms $\calF\circ c$ and $c\circ \calF$ (from 
$\calF^3(3n)$ to $\calF(n)$) coincide}
\end{equation}
\eprop

\cor{ztwo} The categories $\calD^{\mon}(E)$ and 
$\calD^0(E)=\calD^0(\tilE)$ admit a natural action
of the group $\ZZ_2$ (given by the symplectic Fourier-Deligne transform).
\ecor
\ssec{four-on-x}{Fourier-Deligne transforms on the basic affine space}
\sssec{monoid}{Action of a monoid on a category}Let $M$ be a monoid and
let $\calC$ be a category.
\defe{monoid}{\it An action of $M$ on $\calC$} consists of the following
data (\cf. \cite{De1}):

$\bullet$ a functor $F_m:~\calC\to\calC$ for every $m\in M$

$\bullet$ an isomorphism of functors 
$$
\alp_{m_1,m_2}:~F_{m_1}\circ F_{m_2}\to F_{m_1m_2}
$$
for every $m_1,m_2\in M$ such that $m_1m_2$ is defined. 

This data should satisfy the following associativity condition:

\noindent
For every $m_1,m_2,m_3\in M$, such that $m_1m_2m_3$ is defined one
has 
\eq{}
\alp_{m_1m_2}\circ \alp_{m_1,m_2}=\alp_{m_1,m_2m_3}\circ \alp_{m_2,m_3}
\end{equation}
(both sides are isomorphisms between $F_{m_1}\circ F_{m_2}\circ F_{m_3}$
and $F_{m_1m_2m_3}$.
\edefe

We will be interested in the particular case of a {\it braid monoid}
acting on various categories. By definition, a braid monoid associated
to a reductive group $G$ is the semigroup with generators $s_{\alp}$
for every simple root $\alp$ and with braid relations as the only relations.
Braid monoid actions on categories were studied in \cite{De1}.
\sssec{basaff}{The basic affine space of $G$}From now until the end of this
section we assume that $G$ is semisimple and simply connected.
Let $\calB$ denote the flag variety of $G$. By definition this is the variety
of all Borel subgroups in $G$.

We now want to define certain $T$-torsor $\eta: X\to \calB$. 
Unlike the flag variety, the basic affine space $X$  cannot be attached to
the group $G$ canonically. In order to define it we have to fix the 
following data. Recall that $T$ denotes the abstract Cartan group
of $G$. Let $\ome_1,...\ome_n$ denote the fundamental weights of $G$. Then
we must fix the corresponding line bundles $\calO(\ome_1),...,\calO(\ome_n)$
on $\calB$, which are {\it a priori} defined uniquely up
to a non-canonical isomorphism (fixing of $\calO(\ome_i)$ is the
same as fixing of the corresponding representation of $G$ with 
highest weight $\ome_i$, which is also defined uniquely up
to a non-canonical isomorphism). By abuse of notation, we will
denote by $\calO(-\ome_i)$ also the total space of the line bundle
$\calO(-\ome_i)$. Then we define $X$ to be the complement to
``all zero sections'' in $\calO(-\ome_1)\x _{\calB}\x ...\x _\calB 
\calO(-\ome_n)$. Namely, let ${\widetilde\calO(-\ome_i)}$ denote the
complement to the zero section in $\calO(-\ome_i)$. Then we define
\eq{}
X={\widetilde \calO(-\ome_1)} \underset{\calB}\x ...\underset{\calB}\x
{\widetilde \calO(-\ome_n)}
\end{equation}

It is easy to see that 
\eq{}
\calA=\Gam(X,\calO_X)=\bigoplus\limits_{\lam\in P^+(G)}\Gam(\calB,\calO(\lam))
\end{equation}
with the obvious multiplication there (here $P^+(G)$ denotes the
set of integral dominant weights of $G$). The variety $X$ is an open
subset of the affine variety $\oX=\Spec\calA$. The complement
$\oX\backslash X$ is defined by the ideal 
\eq{}
J=\bigoplus\limits_{\lam\in P^{++}(G)}\Gam(\calB,\calO(\lam))
\end{equation}
where $P^{++}(G)$ denotes the set of integral dominant regular weights of
$G$.

It is easy to see, that if we choose a Borel subgroup $B\subset G$ with
unipotent radical $U$, then $X$ can be identified with $G/U$. We will
often make such a choice in order to simplify the discussion.

\sssec{f-on-x}{Fourier transforms on the basic affine space}
For a simple root $\alpha$ let $P_\alpha\subset G$ be the minimal parabolic
of type $\alpha$ containing $B$. Let $B_\alpha$ be the commutator 
subgroup of $P_\alpha$,
and denote $X_\alpha:=G/B_\alpha$. We have an obvious projection
of homogeneous spaces $\pi _\alpha:
X\rightarrow X_\alpha$. 
It is a fibration with the fiber $B_\alpha /U=\A^2-\{0\}$.

Let $\overline \pi _\alpha :\overline X^\alpha \to X_\alpha$ be the relative
affine completion of the morphism $\pi _\alpha$. (So  $\overline \pi _\alpha$
is the affine morphism corresponding to the sheaf of algebras $\pi_{\alpha *}
(\calO_X)$ on $X_\alpha$.) Then  $\overline \pi _\alpha$ has the structure
of a 2-dimensional vector bundle; $X$ is identified with the complement to 
the zero-section in $\overline X^\alpha$. The $G$-action on $X$ obviously
extends to $\overline X^\alpha$; moreover, it is easy to see that the
determinant of the vector bundle $\overline \pi _\alpha$ admits a canonical
(up to a constant) $G$-invariant trivialization, i.e. 
$\overline \pi _\alpha$ admits unique up to a constant $G$-invariant
fiberwise symplectic form $\ome_{\alp}$.  

\noindent
{\bf Examples.}\quad 1) Let $G=SL(2)$ and let $\alp$ denote the unique simple
root of $G$. Then $X\simeq \AA^2\backslash\{ 0\}$, 
$X_\alp=pt$ and $X^\alp\simeq \AA^2$,
i.e. $X^\alp$ is isomorphic to the defining representation of $SL(2)$. 
It is clear that this vector space admits unique up to a constant symplectic
form, which is automatically $G$-invariant.

2) Let $G=SL(n)$. Let $V$ denote the defining representation of $G$, i.e.
$V$ is an $n$-dimensional vector space over $\fqb$ with the natural
$SL(n)$-action. Let us choose some identification of $\Lam^nV$ with $\fqb$, 
i.e. choose a nonzero element $\eps\in \Lam^nV$.
Then the basic affine space of $G$ may be described as follows.
\eq{}
X=\{ (0\subset V_1\subset V_2\subset ... \subset V_n=V,\eps_i\in V_i/V_{i-1})|
\quad \text{such that\ }\eps_1\wedge\eps_2\wedge ...\wedge \eps_n=\eps\}
\end{equation}
Of course, $\eps_n$ is uniquely determined from $\eps_1,...,\eps_{n-1}$ (so
we could omit the $\eps_n$ in the definition, adding the condition
$\eps_i\neq 0$ for $i=1,...,n-1$).

Let $\alp_k$ denote denote the $k$-th simple root of $SL(n)$ ($k=1,...,n-1$).
Then  
\eq{}
\begin{align*}
X_{\alp_k}=\{(0\subset V_1\subset ...\subset V_{k-1}&\subset V_{k+1}\subset ...
\subset V_n=V,\\ 
&0\neq \eps_i\in V_i/V_{i-1}\|\ \text{for $i\neq k,k+1$, and $1\neq i\neq n$} 
\}
\end{align*}
\end{equation}
where $\dim V_i=i$. If we are given some point of $X_{\alp_k}$ as above,
then the fiber of $|pi_{\alp_k}$ over this point is 
$V_{k+1}/V_{k-1} - \{0\}$  (since fixing of $\eps_k$ obviously
fixes $\eps_{k+1}$) and the fiber of
$\opi_{\alp_k}$ is $V_{k+1}/V_{k-1}$. Note that we have canonical
element $\ome \in \Lam^2(V_{k+1}/V_{k-1})$
which is defined by the condition
$$
\eps_1\wedge ...\wedge \eps_{k-1}\wedge \ome\wedge \eps_{k+2}\wedge ...\eps_n=
\eps 
$$
Hence we see that every fiber of $\opi_{\alp_{k}}$ is a $2$-dimensional
vector space with a symplectic structure, which is canonical once we choose
$\eps$.

\smallskip

Let us go back now to an arbitrary group $G$. Let $\alp$ be a simple root
of $G$. Then we may define a functor $\calF_{\alp}:~\calD(X)\to \calD(X)$
by 
\eq{}
\calF_{\alp}=j^*_{\alp}\circ\calF\circ j_{\alp!}
\end{equation}
where $j_{\alp}:~X\to \oX^{\alp}$ is the natural embedding and $\calF$ is
the Fourier-Deligne transform associated with the symplectic vector
bundle $\opi_{\alp}:~\oX^{\alp}\to X_{\alp}$.
 
\sssec{kategorija}{The categories  $\calD^{mon}(X)$ and
$\calD^{reg}(X)$} 
The category $\calD^{\mon}(X)$ is defined analogously
to $\calD^{\mon}(E)$ (cf. Remark \refr{monod} after Definition \refd{mon}).
Note, that ``monodromic'' here means monodromic with respect to the action of
$T$ on $X$. 

\defe{} 1) A tame local system $\calL$ on $T$ is called quasi-regular if
$\alp^*\calL$ is not constant for any coroot $\alp:~\GG_m\to T$.

2) An object $A\in \calD^{\mon}(X)$ is called {\it regular} if
it is glued from $(T,\calL)$-equivariant complexes, where $\calL$ is
a quasi-regular local system on $T$.
\edefe

The following result is proved in the next subsection.

\th{braid} The functors $\calF_{s_\alp}$ map $\calD^{\mon}(X)$ to 
$\calD^{\mon}(X)$. Moreover, there is a canonical action of the braid
monoid $Br_W$, corresponding to $W$ on the category $\calD^{\mon}(X)$ which
extends the functors $\calF_{s_\alp}$.
\eth

We will now define a category $\calD^0(X)$ (which is analogous to the category
$\calD^0(E)$) assuming Theorem \reft{braid}. Namely, it follows from
Theorem \reft{braid} that (since $W$ is canonically embedded into its
braid monoid as a set) for any $w\in W$ we have canonical
functor $\calF_w:~\calD^{\mon}(X)\to \calD^{\mon}(X)$. This functor
maybe constructed as follows. Choose a reduced decomposition
$w=s_{\alp_1}...s_{\alp_k}$ of $w$. Then
$\calF_w=\calF_{s_{\alp_1}}\circ...\circ\calF_{s_{\alp_k}}$. It follows
from Theorem \reft{braid} that $\calF_w$ does not depend on the choice of
a reduced decomposition of $w$ 
up to a canonical isomorphism. We can now define the category $\calD^0(X)$.

\defe{dzerox} The category $\calD^0(X)$ is the full subcategory
of $\calD^{\mon}(X)$, consisting of all the complexes $A\in \calD^0(X)$
such that the following condition holds:

for any simple root $\alp$ of $G$ and for any $w\in W$ the complex
$\calF_w(A)$ lies in $\calD^0(\oX^\alp)$ (recall that $\opi_\alp:~\oX^\alp
\to X_\alp$ is a $2$-dimensional  vector bundle).
\edefe

It is easy to see that one has canonical embeddings
$\calD^{\reg}(X)\subset\calD^0(X)\subset \calD^{\mon}(X)$ of full
subcategories. We will also denote by $\Perv^{\reg}(X)$, $\Perv^0(X)$
and $\Perv^{\mon}(X)$ the corresponding categories of perverse sheaves.
 
One can show (cf. \cite{po}) that the restrictions of the functors
$\calF_w$ to the category $\calD^0(X)$ extend canonically to an action of
$W$ on this category. However, this result is not an immediate corollary
of the above and we will not present the proof in this paper, since
we will not need this in the sequel.
\ssec{r-on-x}{Radon transforms on the basic affine space}

\sssec{}{The space $Z$}Suppose again that we have fixed a basic affine
space $X$ together with the symplectic 
forms $\ome_\alp$ (cf. \refsss{f-on-x}). 
Then we claim that there exists a natural $G$-torsor $Z$ (i.e. 
principal homogeneous space over $G$), endowed with a map
$p:Z\to X$, an action of $T$ and an action of the braid monoid $\brw$, 
corresponding 
to $W$, such that the following
conditions hold:

$\bullet$ $p$ is $G\x T$-invariant 

$\bullet$ the actions of $\brw$ and $T$ on $Z$ agree with the action of 
$W$ on $T$

$\bullet$ the actions of $\brw$ and $G$ on $Z$ commute

First of all let ${Z_{\calB}}$ be the space of all embeddings
of $T$ into $G$. Clearly, this space admits a natural action of
$G\x W$ (where $G$ acts by conjugation on itself and $W$ acts on $T$).
Moreover, it is obvious that the $G$-action is transitive (in fact, a
choice of a point $a:T\hookrightarrow G$ identifies $Z_{\calB}$ with
$G/\Im a$). On the
other hand, we have a natural map $p_{\calB}:{Z_{\calB}}\to \calB$,
defined as follows. Let $a:T\hookrightarrow G$ be an embedding. Then 
there exists a unique Borel subgroup $B$ of $G$, containing the image of $T$,
 such that
the pair $(T,a:T\hookrightarrow B)$ induces the canonical system
of positive roots on $T$ (i.e. $B$ is characterized by the property 
that the set of roots of $T$ on the unipotent radical of the Lie algebra of $B$
is precisely the set of positive roots $R^+$ which the abstract Cartan group
possesses). We now set $p_{\calB}(a)=B$.  

Now we define $Z={Z_{\calB}}\x_{\calB}X$. Then we claim that $Z$
admits a natural action of $G\x \brw$  and an action of
$T$. The actions of $G$ and $T$ on $Z$ are clear: $G$ acts 
on $Z={Z_{\calB}}\x_{\calB}X$ diagonally (since $G$ acts both on
$Z_{\calB}$ and on $X$ this notion makes sense) and $T$ acts just on the
second multiple. The action of $\brw$ on $Z$ may be described for example
as follows.  

First of all, let us describe this action for $G=SL(2)$. In this case
\eq{}
Z=\{x,y\in \AA^2|\ \ome(x,y)=1\}
\end{equation}
and we let the unique non-trivial element $s$ of the Weyl
group of $SL(2)$ act by
\eq{sldva}
s((x,y))=(-y,x)
\end{equation}
(note that $s$ has order $4$). 

For general $G$ we may now do the following.
Assume that we have fixed the fundamental representations $V(\ome_i)$ of
$G$. Then an element $z$ of $Z$ is the same as pair $T\subset B\subset G$
consising of a Borel subgroup $B$ together with a Cartan subgroup
$T$ which is contained in $B$ plus a non-zero vector $v_i\in V(\ome_i)$
which is a highest weight vector of $G$ with respect to $B$ for
every $i=1,...,\rank(G)$. Let now $\alp$ be a simple root of $G$. 
We would like to define 
$s_{\alp}(z)$ (here $s_\alp\in W$ is the corresponding simple reflection.
We set 
\eq{}
s_{\alp}(z)=(T,B^{s_\alp}, v_1^{s_\alp}, ..., v_{\rank G}^{s_\alp})
\end{equation}
where

$\bullet$ $B^{s_\alp}$ is the (unique) Borel subgroup of $G$ which contains
$T$ and is in position $s_\alp$ with $B$

$\bullet$ if $(\ome_i,\alp)=0$ then $v_i^{s_\alp}=v_i$. If 
$(\ome_i,\alp)=1$ then $v_i^{\alp}$ is the unique 
highest weight vector in $V{\ome_i}$
with respect to $B^{s_\alp}$ which satisfies the following property. Let
$\pi_{\calB,\alp}:~\calB\to \calB_\alp$ denote the natural
projection from $\calB$ to the corresponding partial flag variety (thus
the fibers of $\pi_{\calB,\alp}$ are projective lines). 
Let $W_{B,\alp}$ denote the space of global sections of the line bundle   
$\calO(\ome_i)$ restricted to $\pi^{-1}_{\calB,\alp}(B)$. Then 
$W_{B,\alp}$ naturally identifies with the inverse image under
${\overline \pi}_\alp$ of the image of $z$ in the space $X_\alp$
(cf. \refsss{f-on-x} for all the notations). Therefore, 
$W_{B,\alp}$ is a 2-dimensional vector space with has a natural
symplectic form $\ome$ (recall, that we have fixed a $G$-invariant
symplectic form on the bundle ${\overline \pi}_\alp$. Let $w, w^{s_\alp}$
be the images of $v_i$ and $v_i^{s_\alp}$ in $W_{B,\alp}$ respectively.
Then $V^{s_\alp}_i$ is uniquely determined by the requirement 
\eq{}
\ome(w,w^{s_\alp})=1
\end{equation}

It is easy to see that $Z$ satisfies the 
three conditions, listed above. In particular, the $G$-action on $Z$
is simply transitive (this follows immediately from the fact that 
$Z_{\calB}$ is a homogeneous space over $G$, where stabilizer of a point
is a Cartan subgroup, and, on the other hand, $X$ is a $T$-torsor over 
$\calB$).

\sssec{Z-w}{The spaces $Z_w$}
For any $w\in W$ we let $\Gam_w$ denote the graph $w$ on $Z$.
We denote now by $Z_w$ the image of $\Gam_w$ in $X\x X$ (under
the map $p\x p$). 

\prop{rad-everything} 
1) $Z_w$ is a smooth closed subvariety of $X\x X$ of dimension
equal to $\dim X+l(w)$.

2) The restrictions of the projections $\pr_1,\pr_2$ (going from
$X\x X$ to $X$) to $Z_w$ are locally trivial (in \'etale topology) 
fibrations with fiber isomorphic to $\AA^{l(w)}$. 

3) $Z_e$ is equal to the diagonal in $X\x X$ and $Z_{w_0}$ is
isomorphic to $Z$, where the two projections from $Z$ to $X$ are
given by $p$ and $p\circ w_0$ (here $w_0$ denotes the longest 
element in $W$). Also, for any simple reflection $s_\alp\in W$ we
have $Z_{s_\alp}=Z_{{\oX}^{\alp}}$ (recall (cf. \refsss{f-on-x}) that
${\overline \pi_{\alp}}:{\oX}^{\alp}\to X_{\alp}$ is a 2-dimensional
symplectic vector bundle, whose complement to the zero section is
naturally identified with $X$. Recall also that in \refsss{rad-bundle}
we have defined certain variety $Z_E$ (the kernel of the Radon 
transformation) for any vector bundle $E$).  

4) Suppose that $G$ is defined over $\fq$ and let $\fr:X\to X$ denote
the corresponding Frobenius morphism (any $\fq$-rational structure on
$G$ gives rise to a canonical $\fq$-rational structure on $X$). 
Let $\Gam_{\fr}$ denote the graph of the Frobenius morphism on $X$. 
Then the intersection of $Z_w$ with $\Gam_{\fr}$ is transverse. 

5) Let $w_1, w_2\in W$ be such that $l(w_1w_2)=l(w_1)+l(w_2)$. 
Then $Z_{w_1}\circ Z_{w_2}=Z_{w_1w_2}$, where ``$\circ$'' denotes 
composition of correspondences. 
\eprop

The proof of the Proposition is straightforward and it is left to the reader.

\sssec{}{Radon transforms}Let $\pr_{1,2}:X\x X\to X$ denote again the 
two natural projections and let 
$K_w={\overline \QQ}_{\ell,Z_w}[l(w)]\in \calD(X\x X)$.
Then we define the functor $\calR_w:\calD(X)\to \calD(X)$ by
\eq{}
\calR_w(A)=\pr_{2!}(\pr_1^*A\ten K_w)
\end{equation}

\prop{radon-braid} 
1) Let $w_1,w_2\in W$ be such that $l(w_1w_2)=l(w_1)+l(w_2)$. 
Then there is a canonical isomorphism of functors 
\eq{}
\calR_{w_1}\circ \calR_{w_2}\to \calR_{w_1w_2}
\end{equation}

These isomorphisms extend the functors $\calR_w$ to an action of the
braid monoid corresponding to $W$ on the category $\calD(X)$ in the sense
of Deligne (\cite{De1}).

2) The functors $\calR_w$ map the category $\calD^{\mon}(X)$ to itself.
Moreover, for any simple reflection $s_\alp\in W$ the restrictions of 
the functors $\calR_{s_\alp}$ and $\calF_{s_\alp}$ to the category 
$\calD^{\mon}(X)$ are canonically isomorphic.
\eprop

\prf The proof of the first statement follows merely from point 5 in
Proposition \refsss{Z-w}. The proof of the second statement follows 
from Theorem \reft{rad-four} and from the fact that $Z_{s_\alp}$ is
equal to $Z_{\oX^{\alp}}$
(cf. Proposition \refsss{Z-w}(3)).
\epr

\cor{}Theorem \reft{braid} holds.
\ecor

For any $w\in W$ we will denote
by $\Phi_w:\calD^{0}(X)\to \calD^{0}(X)$ the corresponding functor (here we
use the fact that $W$ has a natural set-theoretical embedding into Br$_W$).
One can show (cf. \cite{po}) that the functors $\Phi_w$ extend
canonically to an action of the group $W$ on the category $\calD^0(X)$. 
However, we will not use this in the sequel.

\ssec{}{Remarks on general reductive groups}So far we were speaking 
only about connected simply connected semisimple groups $G$. Let
us explain how to extend the above construction to the case of 
arbitrary connected reductive group $G$. First of all, suppose that
$G$ is semisimple. Let $G^{sc}\to G$ denote the universal
covering of $G$ and let $Z(G^{sc})$ denote the center of $G^{sc}$,
which is the Galois group of the covering $G^{sc}\to G$. Then 
we can repeat all the constructions and statements discussed above,
replacing everywhere the group $G$ by $G^{sc}$ and the words
sheaf on $G/U$" by the words "$Z(G^{sc})$
-sheaf on $G^{sc}/U$". The case of general reductive group
can be treated similarly.

\sec{char}{Character sheaves}

This section is devoted to some 
preliminary material about character
sheaves, that we are going to need in section \refs{proof}. Below $G$ is
a reductive algebraic group defined over $\fq$.
\ssec{characters}{Maximal tori and their characters}

\sssec{}{Maximal tori and conjugacy classes in the Weyl group}
In this section we let $T$ denote the abstract Cartan group of $G$ with
the corresponding $\fq$-rational structure (cf. \cite{dl}, Section 1.1).
Let $W$ denote the corresponding Weyl group, which is endowed with
the action of the corresponding Frobenius automorphism 
$\fr:\ W\to W$. Let $W^{\#}_F$ denote the set of $\fr$-conjugacy
classes of elements in $W$ (by definition, two elements $w_1$ and $w_2$ are
$\fr$-conjugate if there exists $y\in W$ such that 
$w_1=yw_2\fr(y)^{-1})$.

It is shown in \cite{dl} that there is a natural bijection between
$W^{\#}_F$ and the set of conjugacy classes of maximal tori in $G$ which are 
defined over $\fq$. For any $w\in W$ let $T_w$ denote the corresponding
torus (defined uniquely up to $G(\fq)$-conjugacy). 
As an algebraic group over $\fqb$ the torus $T_w$ is identified (canonically)
with $T$. Let $\frw:\ T\to T$ denote the image of the Frobenius morphism
of $T_w$ under this identification. Thus it follows from the definition 
of $T_w$ (cf. \cite{dl}, Section 1.8) that 
\eq{}
\frw=w\circ \fr
\end{equation} 
\sssec{}{Rigidified local systems}Here we follow the exposition of G.~Laumon
(cf. \cite{laum}). Recall that $T$ denotes the abstract Cartan group
of $G$ with the split $\fq$-structure. Let $\calS(T)$ denote the set
of isomorphisms of couples $(\calL,\iota)$, where $\calL$ is an $\ell$-adic
one-dimensional local system on $T$ and $\iota:\ \qlb{\widetilde \to}
\calL_e$ is a rigidification of $\calL$ (here $\calL_e$ denotes the fiber
of $\calL$ at $e\in T$), such that the following property holds.

{\bf Property.} There exists a natural number $n$ such that the pair
$(\mu_n^*\calL,\mu_n^*\iota)$ is isomorphic to the pair
$({\overline \QQ}_{l,T}, {\bf 1})$, where $\mu_n:~T\to T$ is given
by $\mu_n(t)=t^n$,  ${\overline \QQ}_{l,T}$ is
the constant sheaf on $T$ and ${\bf 1}:\ \qlb \to 
({\overline \QQ}_{l,T})_e$
is the obvious map.

The set $\calS(T)$ has canonical structure of an abelian group, given
by tensor product of $\qlb$-sheaves. One has a non-canonical isomorphism
\eq{}
\calS(T)\simeq X^*(T)\ten (\QQ'/\ZZ)
\end{equation}
where $X^*(T)$ is the group of algebraic characters of $T$ and 
\eq{}
\QQ'=\{ {m\over n}\in \QQ|\ m,n\in \ZZ\ \text{and $n$ is invertible in 
$\fq$}\}
\end{equation}

In the sequel we will write simply $\calL\in \calS(T)$, omitting $\iota$
(if it does not lead to a confusion).

The group $W$ acts naturally on the group $\calS(T)$ (an element
$w\in W$ sends $\calL$ to $(w^{-1})^*\calL$). For any $\calL\in \calS(T)$
we let $W_{\calL}$ denote the stabilizer of $\calL$ in $W$ (it is easy to
see that $W_{\calL}$ is independent of the rigidification).

\defe{regularity} (cf. \cite{laum})

1) An element $\calL\in \calS(T)$ is called regular if $W_{\calL}=1$.

2) An element $\calL\in \calS(T)$ is called quasi-regular if for
any coroot $\alp:\ \GG_m\to T$ of $G$ the $\qlb$-sheaf $\alp^*\calL$
on $\GG_m$ is not constant.  
\edefe

\smallskip

{\bf Example.}\quad Let us give an example of a quasi-regular element
of $\calS(T)$, which is not regular. 
Let $G=SL(2)$. In this case $T$ is one-dimensional.
Let $\tau:\ T\to T$ be given
by $\tau(t)=t^2$. Consider $\tau_*({\overline \QQ}_{\ell,T})$. This sheaf
admits a natural action of $\ZZ_2$ (which is the Galois group of the
covering $\tau$). Let $\calL$ denote the skew-invariants of $\ZZ_2$
on $\tau_*({\overline \QQ}_{\ell,T})$. Then $\calL$ is a non-trivial 
one-dimensional local system on $T$ which admits an obvious rigidification.
It is easy to see that $\calL^{\ten 2}\simeq {\overline \QQ}_{\ell,T}$.
Therefore, $\calL$ defines an element in $\calS(T)$. Since
$\calL$ is clearly non-constant, it follows that $\calL$ is
quasi-regular. However, since $\tau$ is a $W$-equivariant map, it
follows that $w^*\calL\simeq \calL$ for any $w\in W$. Hence $\calL$
is not regular. 

However it follows essentially from \cite{dl}, Theorem 5.13 that
if the center of $G$ is connected then the notions of regularity and
quasi-regularity coincide.
\sssec{}{Local systems and characters}Suppose now that our torus $T$
is endowed with some $\fq$-rational structure. Then the Galois group
Gal$(\fqb/\fq)$ acts naturally on $\calS(T)$. Let 
$\fr:\ \calS(T)\to \calS(T)$ denote the 
corresponding Frobenius morphism. Suppose now that we are given
some $(\calL,\iota)\in \calS(T)^{\fr}$. Then $\calL$ acquires a
canonical structure of a Weil sheaf on $T$. Indeed, the fact that 
$(\calL,\iota)\in \calS(T)^{\fr}$ means that there exists an
isomorphism $\calL\simeq \fr^*\calL$ which commutes with $\iota$ and
it is easy to see that such automorphism is automatically unique.
Thus for any $(\calL,\iota)\in \calS(T)^{\fr}$
we may consider the corresponding trace function $\tr(\calL)$ on $T(\fq)$.

\lem{} $\tr(\calL)$ is a character of the finite group $T(\fq)$.
\elem
\ssec{mat}{Matrix coefficients sheaves}Let now $G$ be an arbitrary algebraic
group over $\fqb$, $X$ -- a $G$-space. Take any $A\in \calD(X\x X)$.
Then we can define the matrix coefficient sheaf $m(A)$ by
\eq{matcoef}
m(A)=\pi_!i^*A
\end{equation}
where $i:G\x X\to X\x X$ is defined by $i(g,x)=(gx,x)$ and
$\pi:G\x X\to G$ is the projection on the first variable.
\ssec{modmat}{Modified matrix coefficients}
Take now $G$ again to be
a reductive group over $\fqb$ and $X$ to
be the basic affine space of $G$ and set $\calB=X/T$ to be
the flag variety of G. Let $\Delta:X\to X\x X$ be
the diagonal embedding. In this section we are going to deal
only with those sheaves on $X\x X$, which are 
$(T\x T, \calL\boxtimes \calL^{-1})$-
equivariant for some one-dimensional local system $\calL$ on $T$.
In this case we are going to modify the definition of matrix coefficients.
Namely, since the actions of $G$ and $T$ on $X$ commute, it follows
that the sheaf $i^*A$ in formula \refe{matcoef} is $T$-equivariant
with respect to the action of $T$ on the second variable in $G\x X$.
Therefore, $i^*A$ is a pull-back of some sheaf $\tilA$ on $G\times \calB$.
Let $\tilpi:G\times \calB\to G$ be the projection on the second
variable. We define the {\it modified matrix coefficient sheaf}
by
\eq{}
\tilm(A)=\tilpi_!\tilA
\end{equation}
\ssec{}{Character sheaves}Let us recall Lusztig's definition of
(some of) the character sheaves. Let $\tilG$ denote the variety
of all pairs $(B,g)$, where

$\bullet$\quad $B$ is a Borel subgroup of $G$

$\bullet$\quad $g\in B$

One has natural maps $\alp:\ \tilG\to T$ and $\pi:\ \tilG\to G$
defined as follows. First of all, we set $\pi(B,g)=g$.
Now, in order to define $\alp$, let us remind that for any
Borel subgroup $B$ of $G$ one has canonical identification
$\mu_B:\ B/U_B\wt{\rightarrow} T$, where $U_B$ denotes the unipotent
radical (in fact, this is how the abstract Cartan group $T$ is defined).
Now we set $\alp(B,g)=\mu_B(g)$. 

Let $\calL\in \calS(T)$. We define 
$\chl=\pi_!\alp^*(\calL)[\dim G]$. 
One knows (cf. \cite{lu-char}, \cite{laum})
that the sheaf $\chl$ is perverse, and it is irreducible if
$\calL$ is regular. We want to rewrite
this sheaf as certain (modified) matrix coefficient. Namely,
let $X$ again be the basic affine space of $G$ and let
$\tilDel\subset X\times X$ be the preimage in $X\times X$ of the diagonal 
in $\calB\times\calB$. Then one has canonical morphism
$f:\ \tilDel\to T$ such that for any $(x,y)\in \tilDel$ one
has $y=f((x,y))x$ (recall that $T$ acts on $X$. 

Set now $\al=f^*(\calL)$. 
\lem{matrix-character}
One has canonical isomorphism 
\eq{}
\tilm(\al)[\dim G]\simeq \chl
\end{equation}
\elem

The proof is straightforward. However, we will see in the next section
that this essentially trivial lemma helps us define a Weil
structure on $\chl$ without appealing to the fact that it
is the intermediate extension from the set of regular
semisimple elements in $G$. 
\ssec{weil1}{The Weil structure}
Fix now an isomorphism $\calL\simeq \frw^*(\calL)$ (i.e. we suppose that
$\calL\in \calS(T)^{\frw})$. It was observed by
G.~Lusztig in \cite{lu-char} that fixing such an isomorphism endows 
$\chl$ canonically with a Weil structure. Lusztig's definition of this
Weil structure was as follows. 

Let $j:G_{rs}\to G$ denote the open embedding of the variety of 
regular semisimple elements in $G$ into $G$. 

\lem{}
\eq{}
\chl=j_{!*}(\chl|_{G_{rs}})
\end{equation}

Here $j_{!*}$ denotes the Goresky-McPherson (intermediate) extension (cf.
\cite{bgg}).
\elem

(The lemma follows from the fact that the map $\pi$ is small in the 
sense of Goresky and McPherson).

The lemma shows that it is enough to construct the Weil structure only
on the restriction of $\chl$ on $G_{rs}$. The latter now has a particularly
simple form. Namely, let $\tilG_{rs}$ denote the preimage of
$G_{rs}$ under $\pi$ and let $\pi_{rs}$ denote the restriction of
$\pi$ to $\tilG_{rs}$. Then it is easy to see that 
$\pi_{rs}:\tilG_{rs}\to G_{rs}$ is an unramified  Galois covering with Galois
group $W$. In particular, $W$ acts on $\tilG_{rs}$ and this action is
compatible with the action of $W$ on $T$ in the sense that the restriction
of $\alp$ on $\tilG_{rs}$ is $W$-equivariant.

Now, an isomorphism $\calL\simeq\frw^*(\calL)$ gives rise to an isomorphism
\eq{weil-glupost1}
\alp^*\calL\simeq(w\circ\fr)^*(\alp^*\calL)
\end{equation} 
(here both $w$ and $\fr$
are considered on the variety $\tilG$). Since $\pi_{rs}$ is a Galois covering
with Galois group $W$, it follows that one has canonical identification
\eq{weil-glupost2}
\pi_{rs!}(\fr^*(\alp^*\calL))\simeq  \pi_{rs!}((w\circ\fr)^*(\alp^*\calL))
\end{equation}

Hence from \refe{weil-glupost1} and \refe{weil-glupost2} we get 
the identifications
\eq{}
\fr^*\pi_{rs!}(\alp^*\calL)\simeq \pi_{rs!}(\fr^*(\alp^*\calL))\simeq
\pi_{rs!}((w\circ\fr)^*(\alp^*\calL))\simeq \pi_{rs!}(\alp^*\calL)
\end{equation}
which gives us a Weil structure on $\pi_{rs!}(\alp^*\calL)\simeq 
\chl|_{G_{rs}}$. Hence we have defined a canonical Weil structure on
$\chl$.
\ssec{weil2}{The Weil structure revisited (the case of quasi-regular $\calL$)}
We now want to give a different 
construction of (the same) Weil structure on $\chl$ using lemma 
\refl{matrix-character} and the functors $\Phi_w$. For simplicity we 
will assume here that $\calL$ is quasi-regular, since this is
the only case that we are going to use in the sequel (but, with slight
modifications, the argument presented below works for any local
system $\calL$).

\prop{isow} Let $\calL\in \calS(T)$ be quasi-regular. 
Then one has canonical isomorphism
\eq{isomorphism}
\Phi_{w,w}f^*(\calL)\simeq f^*(w(\calL))
\end{equation}
(recall that $w(\calL)=(w^{-1})^*\calL$).
\eprop

\prf Since the notion of quasi-regularity is invariant under $W$,
it is enough to show that \refe{isomorphism} holds when $w$ is
a simple reflection. In this case, arguing in a standard way we
may assume that $G=SL(2)$.
So, we must construct the isomorphism \refe{isomorphism} when
$G=SL(2)$, $w=s$ -- the unique non-trivial element in the Weyl
group of $SL(2)$ and $\calL$ -- any non-constant element from 
$\calS(T)=\calS(\GG_m)$. 

First of all, let us show that the restriction of $\Phi_{w,w}f^*(\calL)$
to the complement of $\tilDel\subset X\x X$ is equal to zero. Let
$(x,y)\in (X\x X)$. Then the fiber of
$\Phi_{w,w}f^*(\calL)$ at $(x,y)$ can be computed in the following way.

Recall that $X\simeq \AA^2\backslash \{0\}$ and let $Z$ denote the kernel of
the Radon transform on $X$ (with respect to a fixed symplectic form
$\ome$ -- cf. Section \refs{fourier}). 
Let $Z_{x,y}$ denote the closed subvariety of $X\x X$, consisting
of all pairs $(x_1,y_1)\in X\x X$ such that

$\bullet$\quad  $(x_1,y_1)\in \tilDel$

$\bullet$\quad  $(x_1,x)\in Z$   

$\bullet$\quad  $(y_1,y)\in Z$   

Let $j:\ Z_{x,y}\to \tilDel$ be the embedding of $Z_{x,y}$ into $\tilDel$
(or $X\x X$). Then the fiber of $\Phi_{w,w}f^*(\calL)$ at $(x,y)$
is naturally isomorphic to $H^*_c(Z_{x,y},j^*(f^*(\calL))[1]$. Let $H_{x,y}$
denote the stabilizer of $(x,y)$ in $G\x T\x T$.

Suppose now that $(x,y)\not\in \tilDel$. Then the restriction
of $f$ to $Z_{x,y}$ is an isomorphism of the latter variety with $T=\GG_m$.
Indeed, without loss of generality
we may suppose that $\ome(x,y)=1$. Then if $f((x_1,y_1))=\lam$,
i.e. $y_1=\lam x_1$, then the pair $(x_1,y_1)$ is given by
the formulas
$$
x_1=\lam^{-1}x+y, \quad y_1=x+\lam^{-1}y
$$
Hence the pair $Z_{x,y}, j^*(f^*(\calL))$ is isomorphic to the pair
$\GG_m,\calL$ and therefore, since $\calL$ is non-constant, we have
\eq{}
H^*_c(Z_{x,y},j^*(f^*(\calL))=H^*_c(\GG_m,\calL)=0
\end{equation}
Hence the restriction of $\Phi_{w,w}f^*(\calL)$ to the complement of 
$\tilDel$ is equal to $0$.

On the other hand, the restriction of $\Phi_{w,w}f^*(\calL)$ to
$\tilDel$ is equal to $p_!p_1^*(f^*(\calL))[1]$ where we have the
diagram
\begin{center}
   \begin{picture}(7,4)
        \put(0,0){\makebox(1,1){$X\x X$}}
        \put(3,3){\makebox(1,1){$W$}}
        \put(6,0){\makebox(1,1){$X\x X$}}

        \put(3,3){\vector(-1,-1){2}}
        \put(4,3){\vector(1,-1){2}}

        \put(1,2){\makebox(0.5,0.5){$p_1$}}
        \put(5,2){\makebox(0.5,0.5){$p_2$}}
   \end{picture}
\end{center}

Here
$W$ is the subvariety of $X^4$ consisting of all quadruples
$(x_1,x_2,x,y)$, subject to the three conditions above and such
that $(x,y)\in \tilDel$ and
$$
p_1((x_1,y_1,x,y))=(x_1,y_1), \quad p((x_,y_1,x,y))=(x,y)
$$

It is straightforward that $f\circ p_1=s\circ f\circ p$
(recall that $s$ acts by $s(\lam)=\lam^{-1}$ and that $p$
is a smooth fibration with fiber $\AA^1$. This implies
immediately that
\eq{}
p_!p_1^*(f^*(\calL))[1]=f^*(s^*\calL)
\end{equation}
which finishes the proof.
\epr

Let us now explain how the above proposition helps us define a Weil
structure on $\chl$ (for quasi-regular $\calL$). 
For this we need another auxiliary result.

\prop{unitarity}
Let $\calL$ be as above and let $A\in \calD^{\reg}(X\x X)$. Then for any
$w\in W$ one has canonical isomorphisms 
\eq{wmat}
m(A)\simeq m(\Phi_{w,w}A)
\end{equation}

If $A$ is $(T\x T,\calL\boxtimes \calL^{-1})$-equivariant for some
quasi-regular $\calL\in\calS(T)$ then one also has canonical isomorphism
\eq{wmodmat}
\tilm(A)\simeq \tilm(\Phi_{w,w}A)
\end{equation}
\eprop

\prf We will prove here that \refe{wmat} holds. The proof of \refe{wmodmat}
is completely analogous. 

For the sake of simplicity, let us construct a natural isomorphism
of the fibers of $m(A)$ and $\Phi_{w,w}(A)$ at every point $g\in G$. 
Since the functor $\Phi_{w,w}$ ``commutes'' with the $G\x G$-action on
$X\x X$, it
is enough, in fact, to construct such an isomorphism for $g=e$ 
(the unit element), because the fiber of $m(A)$ at the point
$g$ is canonically isomorphic to the fiber of
$m((g\x e)^*A)$ at the point $e$. Let $\Del:~X\to X\x X$ denote the
diagonal embedding of $X$ and let
$\Del X$ denote its image. Then we have to construct a canonical
isomorphism between $H^*_c(\Del^*A)$ and $H^*_c(\Del^*\Phi_{w,w}(A))$. 

Also, 
we may assume, without loss of generality that
$w$ is a simple reflection $s_\alp$ in $W$. It is easy to see that in this 
case our statement reduces immediately to the case $G=SL(2)$.

Let us now make an explicit calculation in this case.
Recall that $X$ in this case is equal to $\AA^2\backslash\{0\}$, which is
endowed with a symplectic form $\ome$,
and we have also the variety $Z$, defined as
$$
Z=\{(x,y)\in X\x X|\ \ome(x,y)=1\}
$$

together with two natural projections $p_1,p_2:~Z\to X$. Let
now $A\in \calD^{\reg}(X\x X)$. Then
$\Phi_{w,w}(A)=(p_2\x p_2)_!(p_1\x p_1)^*A[2]$. Hence,
\eq{}
H^*_c(\Del^*\Phi_{w,w}(A))=H^*_c((p_2\x p_2)^{-1}(\Del X),(p_1\x p_1)^*A)[2]
\end{equation}
(here we denoted by the same symbol $p_1\x p_1$ both
the natural projection from $Z\x Z$ to $X\x X$ and its
restriction to $(p_2\x p_2)^{-1}(\Del X)$).  

Now, the variety $(p_2\x p_2)^{-1}(\Del X)$ can be
described as the set of all triples $(x,y,z)\in X^3$, such
that $\ome(x,z)=1$ and $\ome(y,z)=1$ (the corresponding
point in $Z\x Z$ is $((x,z),(y,z))$). Let now $q$ denote the
restriction of $p_1\x p_1$ to $(p_2\x p_2)^{-1}(\Del X)$. 
Then $q((x,y,z))=(x,y)$. Let us also denote
$(p_2\x p_2)^{-1}(\Del X)$ by $W$. 

Let us compute $H^*_c(q^*A)$ by first applying the functor $q_!$
and then computing cohomology on $X\x X$. Namely, consider
the sheaf $q_!A$. Denote by $j:~Y\to X\x X$  
the embedding of the complement to $\Del X$ into $X$. 
Then we have an exact triangle
\eq{}
j_!j^*(q_!q^*A)\to q_!q^*A\to \Del_{!}\Del^*(q_!q^*A)
\end{equation}

It is easy to see now that $\Del^*(q_!q^*A)$ is
naturally isomorphic to $\Del^*A[-2]$, since the
restriction of $q$ to $q^{-1}(\Del X)$ is a fibration
with fiber $\AA^1$. Therefore, in order to construct an
isomorphism between $H^*_c(\Del^*A)$ and
$H^*_c(\Del^*\Phi_{w,w}(A))=H^*_c(q^*A)[2]$ it is enough
to show that 
\eq{}
H^*_c(j_!j^*(q_!q^*A))=0
\end{equation}
 
Recall now that $\tilDel$ denotes the variety of all pairs
$(x,y)\in X\x X$, which lie on the same line. Let $\tilY$
denote the complement of $\tilDel X$ in $X\x X$ and
let $\tilj:~\tilY\to X\x X$ denote the corresponding embedding.
\lem{}
The canonical map $\tilj_!\tilj^*(q_!q^*A)\to j_!j^*(q_!q^*A)$
is an isomorphism. 
\elem

The lemma follows immediately from the fact that the fiber of
$q$ over any point of $\tilDel X\backslash \Del X$ is empty.

Now we can finish the proof. It follows from the lemma that
it is enough for us to show that $H^*_c(\tilj^*q_!q^*A)=0$. It is
easy to see that $\tilY$ is invariant under the $\GG_m\x \GG_m$-action 
on $X\x X$. However,
we have assumed that $A\in \calD^{\reg}(X\x X)$. On the other hand,
since over $\tilY$ the map $q$ is an isomorphism (which can
be easily checked), it follows 
that $ \tilj^*q_!q^*A=\tilj^*A\in \calD^{\reg}(\tilY)$, i.e. 
$\tilj^*q_!q^*A$ is glued from $(\GG_m\x \GG_m, \calL_1\boxtimes \calL_2)$-
equivariant sheaves, where $\calL_1$ and $\calL_2$ are nontrivial local
systems on $\GG_m$. But $H^*_c(\tilY, B)=0$ 
for any $B\in \calD^{\reg}(\tilY)$,
which finishes the proof.
\epr 

We can now finish constructing the Weil structure on $\chl$. Indeed,
we have
\eq{}
\begin{aligned}
\fr^*\chl\simeq \fr^*\tilm(\al)\simeq \tilm (A(\fr^*\calL))
\simeq\\
\tilm (A(w(\calL)))\simeq \tilm \Phi_{w,w}(\al)\simeq \tilm(\calL)\simeq
\chl
\end{aligned}
\end{equation}
Hence we get a Weil structure on $\chl$.
\ssec{other}{Some other character sheaves}It is easy to see that the 
$G\x T\x T$-orbits 
on $X\x X$
are naturally parametrized by $W$. For $w\in W$ let $\calO_w$ 
denote the corresponding
orbit. 

For any $\calL_1,\calL_2\in \calS(T)$ define
\eq{}
\wll=\{ w\in W|\ w(\calL_1)=\calL_2\}
\end{equation}
Thus
$\wl=W_{\calL,\calL}$. 

In the sequel we will need the
following result.
\prop{other}

1) Let $\calL_1,\calL_2\in \calS(T)$. Then for every $y\in\wll$ 
there is unique (up to isomorphism)
simple 
 $(G\x T\x T,\qlbg\boxtimes\calL_1\boxtimes \calL_2^{-1})$-equivariant
perverse sheaf $\calA_{\calL_1,\calL_2}^y$ 
supported on the closure of $\calO_w$. The sheaves
$\calA_{\calL_1,\calL_2}^y$ form a basis of 
the Grothendieck
group (tensored with $\qlb$) of
the category of $(G\x T\x T,\qlbg\boxtimes\calL_1\boxtimes \calL_2^{-1})$-
equivariant perverse
sheaves on $X\x X$ (or the Grothendieck group of the corresponding
derived category). Therefore, the dimension of this group is equal to
$\#\wll$.

2) If $\calL_1$ and $\calL_2$ are quasi-regular, then for 
every $w_1,w_2\in W$ the sheaf
$\Phi_{w_1,w_2}(\calA_{\calL_1,\calL_2}^y)$ is isomorphic to 
$\calA_{w_1(\calL_1),w_2(\calL_2)}^{w_1^{-1}yw_2}$.

\eprop

\prf For the proof of 1) it is enough to note that $\calO_y$ carries a 
non-zero
$(G\x T\x T,\qlbg\boxtimes \calL_1\boxtimes\calL_2^{-1})$ equivariant 
local system if
and only if $y\in \wll$.

Let us prove 2). Clearly, it is enough to prove
2) when is the pair $(w_1,w_2)$ is of the form 
$(s,e)$ or $(e,s)$ , where $e\in W$ is the unit element and
$s=s_\alp$ is a simple reflection. In this case the statement
of 2) easily reduces to the case $G=SL(2)$ and thus can be checked
by an explicit calculation (as in the proof of Proposition 
\refp{isow}).
\epr

We will denote the sheaf $\calA_{\calL,\calL}^y$ just by $\al^y$. Thus
$\al=\al^e$, where $e\in W$ is the unit element.

\cor{deistvie} Suppose that $\calL\in \calS(T)^{\frw}$ and that $\calL$ is
quasi-regular. Then the perverse sheaf $\Phi_{w,w}(\fr^*\al^y)$ is 
naturally isomorphic
to $\al^{\frw y}$.
\ecor

The proof is straightforward (one should apply Proposition \refp{other}
for $\calL_1=\calL_2=\calL$ and $w_1=w_2=w$).

In the sequel we will need also the following notation. Let 
$\calL_1,\calL_2\in\calS(T)$ and let $y\in \wll$. Then together with
the sheaf $\calA_{\calL_1,\calL_2}^y$ we will consider also the sheaves 
$\calA_{\calL_1,\calL_2}^{y,!}$ and $\calA_{\calL_1,\calL_2}^{y,*}$,
which are respectively ``$!$'' and ``$*$'' extensions of the corresponding
$(G\x T\x T,\qlbg\boxtimes \calL_1\boxtimes \calL_2^{-1}$-equivariant
local system on $\calO_y$ to the whole of $X\x X$. Thus we have
natural morphisms
\eq{}
\calA_{\calL_1,\calL_2}^{y,!}\to \calA_{\calL_1,\calL_2}^y
\to \calA_{\calL_1,\calL_2}^{y,*}
\end{equation}

It is easy to deduce from \cite{lu-char} 
that all character sheaves are by
definition, those perverse sheaves which occur as constituents of sheaves of
the form
$\tilm(A)$, where $A$ is 
$(G\x T\x T,\qlbg\boxtimes \calL\boxtimes \calL^{-1})$-
equivariant. However, we will not need this in the sequel.

\sec{kl}{Kazhdan-Laumon representations: statement of the results}
In this section we are going to define our version of Kazhdan-Laumon
representations and state our main results about them. However, first
we want to explain certain general ideas (which are due to
D.~Kazhdan -- cf. \cite{ka}) which hide behind this definition. We will
explain these general ideas in the case when $G$ is a group over
a local field, but then come back to finite fields again.
\ssec{forms}{Forms of principal series: general ideas}
 
In this section $k$ will be a non-archimedian local field.
For any algebraic variety $Y$ over $k$ we will denote by
$Y(k)$ its set of $k$-rational points. We will assume that our 
$G$ is a semisimple
connected, simply connected split group over $k$. It is an old
ideology (going back to I.~M.~Gelfand) that to every maximal torus
$T$ in $G$, defined over $k$, there should correspond certain ``series'' 
of representations of $G(k)$, parameterized by the characters of $T(k)$.
Existence of such series is also predicted 
by the Langlands local reciprocity law. 
What do we mean by series? In fact, what one would like to 
construct is some canonical smooth representation of the
group $G(k)\times T(k)$ (and then, given a character of $T(k)$, we can
construct a representation of $G(k)$ taking the coinvariants of $T(k)$ with
respect to this character). We will denote this (desired)
representation by $V_{T,k}$.

There is now one case when one (almost) knows the answer for $V_{T,k}$.
This is the case when $T$ is split over $k$. In this case we can take
as $V_{T,k}$ just the corresponding space 
of principal series representations.
I.e. let $X=G/U$ be the basic affine space of $G$, where $U$ is a maximal
unipotent subgroup of $G$ defined over $k$. Then $X$ is a quasi-affine
algebraic variety, defined over $k$, and we take $V_{T,k}=\calS(X(k))$ --
the space of locally constant compactly supported functions on the
set $X(k)$ of $k$--rational points of $X$ (note that $X(k)=G(k)/U(k)$).
Since $G(k)$ acts naturally (on the left) on $X(k)$ and $T(k)$ acts
there naturally on the right ($T$ is a torus in $G$ which normalizes $U$),
our space is a representation of $G(k)\x T(k)$.

What can we do for other tori? Let us describe one idea in this direction
(cf. \cite{ka}). Galois theory 
tells us that different conjugacy classes of maximal tori in $G(k)$ are
classified by the homomorphisms $\Gam\to W$ (up to $W$-conjugacy)
where $\Gam$ is the absolute
Galois group of $k$ and $W$ is the Weyl group of $G$. Now one would like
to think about different series of representations of $G(k)$ as ``forms''
of the principal series. Let us see how we can do it. Let $T$ be a torus
in $G$ defined over $k$. 
We will assume the existence of $V_{T,k}$ and try to see
what it implies. The Langlands reciprocity law tells us that for any
Galois extension $k'/k$ there should exist a ``lifting'' of $V_{T,k}$
to a representation of $G(k')\x T(k')$ on which the group Gal$(k'/k)$
acts, and this representation should be isomorphic to $V_{T,k'}$.
Suppose now that $T$ splits over $k'$. Then we must get some non-trivial
action of Gal$(k'/k)$ on $\calS(X(k'))=V_{T,k'}$. How to construct it?
Since $T$ splits over $k'$, it corresponds to a homomorphism
$\pi:\text{Gal}(k'/k)\to W$. Therefore, it is enough for us to 
describe an action of the group $W$ on $\calS(X(k'))$ (and then
twist the obvious action of Gal$(k'/k)$ on this space (coming
from its action on $X(k')$) by means of $\pi$). This is already
an ``algebraic'' problem (i.e. it has nothing to do with field
extensions), and so, it is enough to describe the $W$-action on 
the space $\calS(X(k))$. Unfortunately, the space $\calS(X(k))$
does not admit any natural action of $W$. However, it is well-known
that the space $L^2(X(k),\mu)$ does, where $\mu$ is a $G(k)$-invariant
measure on $X(k)$. Consider, for simplicity, the example $G=SL(2)$.
Then $X=\AA^2\backslash\{0\}$ and $L^2(X(k))=L^2(k^2)$. Also 
$W=\pm 1$. The space $k^2$ admits unique up to a scalar 
$SL(2,k)$-invariant symplectic form $\ome$. Let $-1\in W$ act
by the Fourier transform $F$ on $L^2(k^2)$ where we identify  $k^2$
with the dual vector space by means of $\ome$. Then $F^2=id$ (because
$\ome$ is symplectic) and this is the desired action. In the general
case the corresponding action can be described just repeating the
construction described in \refss{four-on-x} and replacing everywhere
sheaves by functions 
(cf. \cite{ka}). 

The first question which immediately appears here is the following

\noindent
{\bf Question 1.} Find a ``nice'' $W$-invariant subspace of $L^2(X(k))$
which contains $\calS(X(k))$. This is an important question
in itself. In particular, the definition of the above space should make
sense in the adelic situation as well, where it might be used to
study analytic properties of Eisenstein series. Also this
space should give the ``correct'' $W$-equivariant version of principal
series (this is why I wrote that one ``almost'' knows
the answer for $V_{T,k}$ in the case of a split torus). 
Of course, one could take just the minimal subspace 
with the above properties. But such a definition is very difficult
to work with (in particular, it does not make sense in the global situation).
On the other hand for $G=SL(2)$ this space is just the space $\calS(k^2)$
of smooth compactly supported functions on $k^2$. One would like to
have a similar description of this space for any $G$. We will not discuss
a solution of this question in this paper. Instead, (in the case of finite 
fields) we will restrict ourselves to some smaller space, on which
$W$ will act. This space will be roughly ``the space of all trace
functions of Weil sheaves lying in $\calD^0(X)$'').

But even if we can answer Question 1, then we arrive to a much more
intriguing

\noindent
{\bf Question 2.} How to construct the representations $V_{T,k}$ 
which would satisfy all the above properties? At the moment an
answer is given only in the case when we replace our field $k$
by a finite field (in that case, of course, one has the Deligne-Lusztig
construction of representations, but in view of what was said above
one would like to have a different construction, which is more
``compatible'' with the theory of forms of algebraic objects. A discussion
of such a construction is presented in 2.2). In the $p$-adic case
D.~Kazhdan has given some conjectural way to construct forms of
the principal series for $G=GL(n)$ (however, he did not prove that
these forms are really well-defined -- cf. \cite{ka}.
\ssec{notations}{Notations}Now we come back to the case of finite field.
>From now on in this section we will suppose 
that the group $G$ is defined over $k=\fq$. 
Our purpose here is to define certain representations
$V_{\calL,w}$ of the finite group $\gfq$, where $\calL\in \calS(T)^{\frw}$ 
for some $w\in W$ (moreover we will see that, in fact, $V_{\calL,w}$ depends
only on $\calL$ and not on $w$). The definition will 
be a slight modification of that from
\cite{kl}. We hope that the connection of the definition with the
general discussion of the previous section will be clear.

In the definitions below will work with Grothendieck groups of various
abelian categories of perverse sheaves on $X$, endowed with additional 
structures. However, one can also work with the Grothendieck group of the 
corresponding derived category, which would (trivially) lead
to the same answer. 
\ssec{}{Kazhdan-Laumon representations}
\sssec{w=1}{Some motivation}
Before giving the definition, let us explain some general fact
(which is quite well-known -- cf. \cite{kl}). 
This fact will serve as a motivation for the 
definition. Namely, following the ideas, described in \refss{forms}
we need to define ``forms'' of the space of functions on
$X(\fq)$. For this we need to describe a more algebraic definition
of this space.

Let $X$ be an arbitrary scheme of finite type over $\fq$ and let $L(X(\fq))$
denote the space of functions on the finite set $X(\fq)$. We want to give 
some purely algebro-geometric construction of this finite-dimensional
space, which will not appeal also to the notion of $\fq$-rational point of
$X$. 

Let $K(X)$ denote the Grothendieck group of perverse sheaves on $X$, tensored
with $\qlb$ (which is the same as the Grothendieck group of the corresponding
derived category, tensored with $\qlb$. However, for simplicity, 
we prefer to work with
abelian categories, and not with triangulated ones). We have a natural 
surjective
map
$\tr:K(X)\to L(X(\fq))$. We would like to identify the kernel of this
map. This is done as follows. One can define a canonical $\qlb$-valued
symmetric pairing $\<\cdot,\cdot\>:K(X)\ten K(X)\to \qlb$ in the following
way. Let $(A,\alp:A\to \fr^*A),(B,\beta:B\to \fr^*B)$ be two Weil (perverse)
sheaves on $X$. 
Then we have an induced 
endomorphism $\phi(\alp,\beta):\RHom(A,DB)\to \RHom(A,DB)$ defined as
the composition
\eq{}
\RHom(A,DB)\to \RHom(\fr^*A,\fr^*DB)\to \RHom(A,DB)
\end{equation}
and we define
\eq{}
\<(A,\alp),(B,\beta)\>=\sum\limits_i(-1)^i\Tr(\phi(\alp,\bet),\Ext^i(A,DB))
\end{equation}
It is clear, that$\<\cdot,\cdot\>$ descends to a well-defined pairing
on $K(X)$. The following result can be easily deduced from the 
Grothendieck-Lefschetz trace formula for the Frobenius correspondence.
\prop{}
The kernel of $\tr:K(X)\to L(X(\fq))$ is equal to the kernel of
the pairing $\<\cdot,\cdot\>$.
\eprop

Hence we obtain a canonical isomorphism between $L(X(\fq))$ and 
$K(X)/\Ker\<\cdot,\cdot\>$.

\sssec{w_arb}{Definition of Kazhdan-Laumon representations}Let now 
$\calL\in\calS(T)^{\frw}$, i.e. 
$\calL$ is a local system on $T$ such that $\frw^*\calL$ is isomorphic
to $\calL$. Our goal is to define certain $\ell$-adic representation
$V_{\calL,w}$ of the group $\gfq$ (we will see afterwards that $V_{\calL,w}$
will actually depend only on $\calL$ and not on $w$). We will use the
notations of \refss{four-on-x}. 

Let $\Perv^0_{\calL,w}(X)$ denote the category, whose objects are pairs
$(A,\alp)$, where 

$\bullet$\quad $A$ is a $(T,\calL)$-equivariant perverse sheaf on $X$, which
lies in $\Perv^0(X)$ (as an abstract perverse sheaf)

$\bullet$\quad $\alp:A\simeq \Phi_w(\fr^*A)$ is an isomorphism

(note that an isomorphism $\calL\simeq \fr^*\calL$ endows 
$\Phi_w(\fr^*A)$ with the structure of $(T,\calL)$-equivariant sheaf).

Morphisms in the category $\Perv^0_{\calL,w}(X)$ are morphisms between sheaves,
which commute with $w$. Note that if $(A,\alp)\in\Perv^0_{\calL,w}(X)$ then
$(DA,D(\alp)^{-1})\in \Perv^0_{\calL^{-1},w}(X)$ (here we use the fact that
$\Phi_w$ commutes with Verdier duality). The 
following lemma is easy (it follows from the 
exactness of the functor $\Phi_w$ on the category $\Perv^0(X)$). 
\lem{} The category $\Perv^0_{\calL,w}(X)$ is abelian.
\elem
 
It is easy to see that the category $\Perv^0_{\calL,w}(X)$ admits
a natural action of the group $\gfq$ (coming from the geometric
action of this group on $X$).

Set now $K_{\calL,w}=K(\Perv^0_{\calL,w}(X))\ten \qlb$ (here
$K(\Perv^0_{\calL,w}(X))$ is the Grothendieck group of the category
$\Perv^0_{\calL,w}(X)$. 
This is an infinite-dimensional vector space over $\qlb$, on which the finite
group $\gfq$ acts. 
We now want to define certain quotient of this $\gfq$-representation, which 
will already be finite-dimensional. 

First of all, we claim that there is a natural $\gfq$-invariant
pairing $\la\cdot,\cdot\ra$ between $K_{\calL,w}$ and $K_{\calL^{-1},w}$.
It is constructed in the following way. 
Let $(A,\alp)\in \Perv^0_{\calL,w}(X)$ and 
$(B,\bet)\in\Perv^0_{\calL^{-1},w}(X)$. Consider $\RHom(A,DB)$ (the $\RHom$
is computed just in the category $\calD(X)$). Then we have an induced 
endomorphism $\phi(\alp,\beta):\RHom(A,DB)\to \RHom(A,DB)$ defined as
the composition
\eq{}
\RHom(A,DB)\to \RHom(\Phi_w(A),\Phi_w(DB))\to \RHom(A,DB)
\end{equation}

\rem{} Note that in section \refs{trace} we will use 
$\phi(\alp,\beta)$ to denote a slightly different endomorphism.
\erem

Set now
\eq{}
\Tr(\phi(\alp,\bet))=\sum\limits_i (-1)^i\Tr(\phi(\alp,\bet),\Ext^i(A,DB))
\end{equation}
We now define $\la (A,\alp),(B,\bet)\ra= 
\frac{\Tr(\phi(\alp,\bet))}{\# T_w(\fq)}$. It is clear
that $\la\cdot,\cdot\ra$ descends to a well-defined pairing between 
$K_{\calL,w}$ and $K_{\calL^{-1},w}$, which is $\gfq$-equivariant.
Denote now by $K^{\text{null}}_{\calL,w}$ the left kernel of the pairing
$\la\cdot,\cdot\ra$, i.e. 
\eq{}
K^{\text{null}}_{\calL,w}=\{ a\in K_{\calL,w}|\ 
\la a,b\ra=0\ \text{for any}\ b\in K_{\calL^{-1},w}\}
\end{equation}
Set now $V_{\calL,w}=K_{\calL,w}/K^{\text{null}}_{\calL,w}$.
It is clear that the group $\gfq$ acts on $V_{\calL,w}$. However, 
it is not clear {\it a priori} whether this representation of 
$\gfq$ possesses any good properties (for example, it is not obvious that
$V_{\calL,w}$ is finite-dimensional).
The following theorem says, in a sense, that it is really the case.
It is one of the main results of this paper. 
\th{main}

1) $V_{\calL}$ is finite-dimensional. Moreover, if $\calL$ is quasi-regular,
then 
\eq{dimhom}
\Hom_{\gfq}(V_{\calL,w},V_{\calL,w})=\# W_{\calL}^{\frw}
\end{equation}
Hence, if $\calL$ is regular, then
$V_{\calL}$ is irreducible.

2) Suppose that $\calL$ is quasi-regular. Then 
the character of $V_{\calL}$ is equal to $\tr(\chl)$ (cf.
subsection \refss{char-sheaves})  
\eth

\noindent
\rem{} One should compare 2) with \cite{dl}, Theorem 6.8 (where
analogous statement for Deligne-Lusztig representations is proved). In fact,
one can show also the following generalization of \refe{dimhom}:
\eq{generalization}
\dim\Hom_{\gfq}(V_{calL_1,w_1},V_{\calL_2,w_2})=\#
\{ w\in W|\ w(\calL_2)=\calL_1, \ \text{and}\ w_1^{-1}\fr(w)w_2=w\}
\end{equation}
We will sketch the proof of \refe{generalization} 
in the next section.
\erem

\noindent
{\bf Example.}\quad Let $w=1$ and suppose for simplicity that $\calL$ is 
quasi-regular. Let $\theta=\tr(\calL)$ be the corresponding character
of $T(\fq)$. Then it is easy to see from \refsss{w=1} that in this case
we have
\eq{}
V_{\calL}=\{ \phi:X(\fq)\to\qlb|\ f(xt)=\theta(t)f(x)\ 
\text{for any $x\in X(\fq),
t\in T(\fq)$}\}
\end{equation}
i.e. $V_{\calL}$ is just the corresponding principal series representation.  

\sec{proof}{Proof of Theorem \reft{main}}

\subsection*{} This section is devoted to the proof of Theorem \reft{main}, as
well as to some generalization of part 2 of this theorem.

\ssec{}{Whittaker sheaves}This subsection is devoted to the proof that 
$V_{\calL,w}$ is non-zero. Moreover, we will give an
estimate of $\dim\Hom_{\gfq}(V_{\calL,w},V_{\calL,w})$ from below.
\sssec{}{Whittaker category}
Let us choose an $\fq$-rational maximal unipotent subgroup
$U$ in $G$. For any simple root $\alp$ of $G$ let us denote
by $U_{\alp}$ the corresponding one-parameter subgroup of $U$.
\defe{} A character (homomorphism of algebraic groups over
$\fq$) $\eps:~U\to \GG_a$ is called
{\it non-degenerate} if 
$\eps|_{U_\alp}$ is non-trivial for every simple root $\alp$ of $G$ (here
$U_\alp$ denotes the one-parametric subgroup of $U$, corresponding
to $\alp$).
\edefe

Fix a non-degenerate character $\eps$ of $U$. Let $\calL(\psi,\eps)$
denote the one-dimensional local system $\eps^*\calL_{\psi}$ on $U$ 
(recall that in Section \refs{fourier} we have fixed a non-trivial additive
character $\psi$ of $\fq$ and that we denote by $\calL_{\psi}$ the 
corresponding Artin-Schreier sheaf). Let also
$\calD_{\eps,\psi}(X)$ denote 
the derived category of $(U,\calL(\eps,\psi))$-equivariant
sheaves on $X$. The functors $\Phi_w$ clearly extend to $\calD_{\eps,\psi}(X)$.
\prop{whittaker}

1) The image of $\calD_{\eps,\psi}(X)$ under the forgetful functor
to $\calD(X)$ lies in $\calD^0(X)$.

2)The category $\calD_{\eps,\psi}(X)$ is equivalent to the category
$\calD(T)$. Moreover, for any choice of the symplectic forms
$\ome_i$ (cf. Section \refs{fourier}) there exists a choice of $\eps$ and
the above equivalnece of categories, in such a way that under 
this equivalence the functors $\Phi_w$
will tranform into the geometric action of $W$ on $\calD(T)$. Moreover,
this equivalence commutes with the functor $\fr^*$. 
\eprop

This statement is proven in \cite{ka} on the level of functions.
Here one should just repeat word-by-word the arguments of \cite{ka}.

\sssec{}{Whittaker model of Kazhdan-Laumon representations}
In what follows we fix $\eps$ and an equivalence 
$\calD_{\eps,\psi}(X)\simeq \calD(T)$ which satisfy the conditions
of Proposition 
\refp{whittaker}. Let $B$ denote the Borel subgroup, ccontaining $U$. Let
$j:~C\to X$ denote the embedding of the unique open
$B$-orbit on $X$. $C$ has a natural action of $B\x T$ and
it can be $U\x T$-equivariantly identified with $U\x T$.
Choose now $\calL\in \calS(T)$ and set $W_{\calL,\eps,\psi}=
j_!(\calL(\eps,\psi)\boxtimes \calL)$. The following lemma
is straighforward.

\lem{}$W_{\calL,\eps,\psi}$ is the unique (up to an isomorphism)
irreducible 
$(U\x T,\calL(\eps,\psi)\boxtimes \calL)$-
equivariant perverse sheaf on $X$. 
\elem

It follows from Proposition \refp{whittaker} that if 
$\calL\in \calS(T)^{\frw}$ then there is a natural isomorphism 
\eq{}
\alp_{\eps,\psi,\calL}:~W_{\calL,\eps,\psi}{\widetilde\to}\fr^*
\Phi_w(W_{\calL,\eps,\psi})
\end{equation}

and also $(W_{\calL,\eps,\psi},\alp_{\eps,\psi,\calL})\in \Perv^0_{\calL,w}$.
We claim now that the image of $(W_{\calL,\eps,\psi},\alp_{\eps,\psi})$ in
$V_{\calL,w}$ is non-zero. Indeed, it is easy to see that
the canonical pairing $\<,\>$ of $W_{\calL,\eps,\psi}$ with
its dual is equal to $1\neq 0$ which proves
what we want. 

Thus $V_{\calL,w}\neq 0$. Working more accurately (varying $\eps$)
we can show also that
\eq{whit-neravenstvo}
\dim\Hom_{\gfq}(V_{\calL,w},V_{\calL,w})\geq \#\wl^{\frw}
\end{equation}

\ssec{}{Proof of theorem \reft{main}(1)}
In this subsection we fix $\calL$ and $w$ and we will write
$V_{\calL}$ instead of $V_{\calL,w}$.
Let us show that $V_{\calL}$ is 
finite-dimensional for any $\calL\in \calS(T)$. Since we have a
canonical perfect $\gfq$-invariant pairing $\<\cdot,\cdot\>$ between
$V_{\calL}$ and $V_{\calL^{-1}}$, it is enough to show that
\eq{mainestimate}
\dim (V_{\calL^{-1}}\ten V_{\calL})^{\gfq}\leq \#\wl
\end{equation}

This is done in the following way.
\lem{subquotient} One can identify canonically the $\gfq\x \gfq$-module 
$V_{\calL^{-1}}\ten V_{\calL}$ with a subquotient of 
$V_{\calL^{-1}\boxtimes \calL}$. Here $V_{\calL^{-1}\boxtimes \calL}$ is
the representation of $\gfq\x \gfq$ constructed as in the previous section,
using the local
system $\calL^{-1}\boxtimes \calL$ on $T\x T$ (note that the basic affine
space of $G\x G$ is $X\x X$).
\elem
\prf 
First of all, there is a natural map 
$K_{\calL^{-1},w}\ten K_{\calL,w}\to K_{\calL^{-1}\boxtimes \calL}$
constructed as $(A,\alp)\ten (B,\beta)\mapsto 
(A\boxtimes B,\alp\boxtimes \bet)$. It follows easily from the
K\"uneth formula that this map preserves the natural pairings on
$K_{\calL^{-1},w}\ten K_{\calL,w}$ and $K_{\calL^{-1}\boxtimes \calL}$.
Hence the kernel of the composite map
$K_{\calL^{-1},w}\ten K_{\calL,w}\to V_{\calL^{-1}\boxtimes\calL}$ lies in 
$K^{\text{null}}_{\calL^{-1}}\ten K_{\calL}+
K_{\calL^{-1}}\ten K^{\text{null}}_{\calL}$. Hence, 
$V_{\calL^{-1}}\ten V_{\calL}$ gets identified with a subquotient of
$V_{\calL^{-1}\boxtimes \calL}$.

\epr

The lemma shows that it is enough to show that 
$\dim V_{\calL^{-1}\boxtimes \calL}^{\gfq}\leq \# \wl^{\fr}$. 
\lem{average}
Any element of $V_{\calL^{-1}\boxtimes \calL}^{\gfq}$ can be represented
by a $G$-equivariant pair $(A,\alp)$. 
\elem

\prf
Indeed, let $v\in V_{\calL^{-1}\boxtimes \calL}$ be
represented by some $(\oA,\oalp)\in\calD^0_{\calL,w}(X)$. For any $G$-variety
$Y$ let 
$\calD_G(Y)$ denote the derived category of equivariant constructible
$\qlb$-sheaves on $Y$ (cf. \cite{bl}). 
Let also $\avg:\ \calD(Y)\to\calD_G(Y)$ denote
the functor of ``averaging'' over $G$. By definition, $\avg(A)=p_!a^*A$
where the maps $a$ and $p$ are defined by the following commutative diagram:
\begin{center}
   \begin{picture}(7,4)
	\put(0,0){\makebox(1,1){$Y$}}
	\put(3,3){\makebox(1,1){$G\x Y$}}
	\put(6,0){\makebox(1,1){$Y$}}

	\put(3,3){\vector(-1,-1){2}}
	\put(4,3){\vector(1,-1){2}}

	\put(1.5,2){\makebox(0.5,0.5){$a$}}
	\put(5,2){\makebox(0.5,0.5){$p$}}
   \end{picture}
\end{center}

Here $p$ is the projection on the second variable and $a$ is given by
$a((y,g))=g^{-1}y$. The functor $\avg[2\dim G]$ is left 
adjoint to the forgetful functor $\calD_G(Y)\to \calD(Y)$.

Let us now go back to the case when $Y=X\x X$. We set now
\eq{}
(A,\alp)= (\avg(\oA),\,  {\avg(\oalp)\over \#\gfq})
\end{equation}
(note that $\avg(\oalp)$ is defined, since $\Phi_{w,w}$ commutes with 
the $G$-action).
 
It follows now easily from the usual Grothendieck-Lefschetz trace 
formula (for the Frobenius morphism) that the image of $(A,\alp)$ 
in $V_{\calL^{-1}\boxtimes \calL}$ is equal to that of 
$(\oA,\oalp)$, which is $v$. On
the other hand, by definition, the pair $(A,\alp)$ is $G$-equivariant.
This finishes the proof. 
\epr

\sssec{}{End of the proof}We claim now that Proposition \refp{other}
together with Corollary \refc{deistvie} imply our statement. Indeed,
let $y\in \wl^{\frw}$. Then the sheaf $\al^y$ defines us a $G$-equivariant
object
in $\Perv^0_{\calL^{-1}\boxtimes \calL,(w,w)}(X\x X)$ (by Corollary
\refc{deistvie}).
Hence it defines an element $a_{\calL,w}^y$ in 
$(V_{\calL,w}\otimes V_{\calL^{-1},w})^{\gfq}$
and it follows from Proposition \refp{other} and from Lemma \refl{average}
that the elments $a_{\calL,w}^y$ (for all $y\in \wl$) span 
$(V_{\calL,w}\otimes V_{\calL^{-1},w})^{\gfq}$. 
Therefore we see that 
\eq{neravenstvo}
\dim V_{\calL^{-1}\boxtimes \calL}^{\gfq}\leq \#\wl^{\frw}
\end{equation}

which by what is explained above implies that 
$V_{\calL}$ is finite-dimensional and that 
\eq{ravenstvo}
\dim\Hom(V_\calL,V_\calL)=\#\wl^{\frw}
\end{equation}
(the last equation follows from \refe{neravenstvo} and from \refe{whit-neravenstvo}).

We claim now that a little more is true. Namely, we claim that the subquotient
which appears in the formulation of Lemma \refl{subquotient} coincides with 
the whole of $V_{\calL^{-1}\boxtimes\calL}$. Since both 
$V_{\calL^{-1}}\ten V_{\calL}$ 
and
$V_{\calL^{-1}\boxtimes \calL}$ are finite-dimensional
representations of the group  $\gfq\x \gfq$ (the fact that $V_{\calL^{-1}\boxtimes \calL}$
is finite-dimensional
follows from the same arguments as above, applied to the group $G\x G$), it is enough
to show that
\eq{idiotism}
\dim\Hom_{\gfq\x\gfq}(V_{\calL^{-1}}\ten V_{\calL},V_{\calL^{-1}}\ten V_{\calL})=
\dim\Hom_{\gfq\x\gfq}(V_{\calL^{-1}\boxtimes \calL},V_{\calL^{-1}\boxtimes \calL})
\end{equation}

However the left hand side of \refe{idiotism} is equal to
$(\#\wl^{\frw})^2$ by \refe{ravenstvo} and the right hand side of \refe{idiotism} is
equal to $(\#\wl^{\frw})^2$ again by \refe{ravenstvo}, but applied to the 
group $G\x G$ and the $T\x T$ local system $\calL\boxtimes \calL^{-1}$. 
Hence \refe{idiotism}
holds, which finishes the proof.

We will denote by $\theta:~V_{\calL^{-1}}\ten V_{\calL}{\widetilde\to}V_{\calL^{-1}
\boxtimes \calL}$
the resulting isomorphism.
\ssec{}{The inverse map}We now want to construct an inverse map
$\phi:~V_{\calL^{-1}\boxtimes \calL}\to V_{\calL^{-1}}\ten
V_{\calL}$.

Let $\calK\in \calD_{{\calL^{-1}\boxtimes \calL},(w,w)}(X\x X)$. Then we define
a functor $\Phi_{\calK}:~\calD_{\calL}\to \calD_{\calL}$ by putting
\eq{}
\Phi_{\calK}(A)=p_{2!}(p_1^*A\ten \calK)
\end{equation}

where $A\in \calD_{\calL}(X)$ and $p_1,p_2:~X\x X\to X$ are the natural
projections. 

Since we have assumed that 
$\calK\in \calD_{{\calL^{-1}\boxtimes \calL},(w,w)}(X\x X)$ it
follows that we are given an isomorphism of functors
\eq{frobenius}
\Phi_{\calK}\simeq (\fr^*\circ\Phi_w)\circ \Phi_{\calK}\circ (\fr^*\circ\Phi_w)
\end{equation}

Hence if we are given $A\in \calD_{\calL,w}(X)$, i.e. $A$ is endowed
with an isomorphism $\fr^*\Phi_w(A)\simeq A$ then using \refe{frobenius}
we may also identify $\Phi_{\calK}(A)$ with 
$\fr^*\Phi_w(\Phi_{\calK}(A))$.  Thus $\Phi_{\calK}$ induces a
map from $K_{\calL,w}$ to $K_{\calL,w}$ which we will denote by
$\phi_{\calK}$. 

\prop{inversemap} 
Let $\calK\in \calD_{\calL^{-1}\boxtimes \calL, (w,w)}(X\x X)$.
Then
\begin{enumerate}
\item $\phi_{\calK}(K^{\text{null}}_{\calL,w})\subset 
K^{\text{null}}_{\calL,w}$
\item Suppose that the image of $\calK$ in 
$K_{\calL^{-1}\boxtimes \calL,(w,w)}$ lies in 
$K^{\null}_{\calL^{-1}\boxtimes \calL,(w,w)}$. Then 
\eq{}
\phi_{\calK}(K_{\calL,w})\subset 
K^{\text{null}}_{\calL,w}
\end{equation}
\end{enumerate}
\eprop

Proposition \refp{inversemap} implies that the assignement
$\calK\to \phi_{\calK}$ descends to a well defined map
$\phi:~V_{\calL^{-1}\boxtimes \calL,(w,w)}\to \End V_{\calL,w}=
V_{\calL^{-1},w}\ten V_{\calL,w}$.

\prf Let us prove (2). The proof of (1) is analogous. 

So, suppose that we are given $\calK$ as above whose image in
$V_{\calL^{-1}\boxtimes \calL}$ vanishes. In particular, 
$\calK$ is endowed with an isomorphism $\gam:\calK{\widetilde \to}
\fr^*\Phi_{w,w}(\calK)$.

Let now $(A,\alp),(B,\beta)\in \calD_{\calL,w}(X)$. Then by \refe{frobenius}
the complex $\Phi_{\calK}(A)$ is also endowed with an isomorphism
$\Phi_{\calK}(A)\simeq \fr^*\Phi_w(\Phi_{\calK}(A)$ which we will denote
by $\Phi_{\calK}(\alp)$. Thus we may consider the endomorphism
$\phi(\Phi_{\calK}(\alp),\beta)$ of $\RHom(\Phi_K(A),B)$. We must show
that $\Tr(\phi(\Phi_{\calK}(\alp),\beta))=0$. 

On the other hand, one has
\eq{}
\begin{align*}
\RHom(\Phi_{\calK}(A)&,B)=\RHom(p_{2!}(p_1^*A\ten \calK),B)=\\
&\RHom(p_1^*A\ten \calK,p_2^!B)=\RHom(\calK, DA\boxtimes B)
\end{align*}
\end{equation}
It is easy to see that under this identifications the map
$\phi(\Phi_{\calK}(\alp),\beta)$ goes to 
$\phi(\gam, \alp\boxtimes \beta)$ (the latter is an endomorphism
of $\RHom(\calK, DA\boxtimes B)$. 

Now, the fact that $\calK$ vanishes in $V_{\calL^{-1}\boxtimes \calL}$
implies that $\Tr(\phi(\gam, \alp\boxtimes \beta))=0$. Hence
$\Tr(\phi(\Phi_{\calK}(\alp),\beta))=0$, which finishes the proof.
\epr

\lem{}
\begin{enumerate}
\item
$\phi\circ\theta =\#T_w(\fq)\id$
\item 
Let $m(\calK)$ denote the matrix coefficient sheaf, corresponding
to $\calK$ (cf. \refss{mat}) with the Weil structure, defined as in \refss{weil2}.
Then the $\qlb$-valued function $\frac{\tr(m(\calK))}{\#T_w(\fq)}$ is equal 
to the matrix coefficient
of the image of $\calK$ in $V_{\calL^{-1}\boxtimes \calL}=V_{\calL^{-1}}\ten V_{\calL}$. 

\end{enumerate}
\elem

The lemma is proved by a direct computation, which is left to the 
reader.

\ssec{}{A variant} Let us again assume that we are given 
$\calK\in \calD_{\calL^{-1}\boxtimes \calL,(w,w)}$. Then, 
analogously to the definition of modified matrix coefficients,
we may define a slightly different version of the functor
$\Phi_{\calK}$ which we will denote by $\tilPhi_{\calK}$.
Namely, suppose that we are given 
$A\in \calD_{\calL}(X)$. Then there exists some $\calF\in \calD(\calB\x X)$
whose pull-back to $X\x X$ can be identified with $p_1^*A\ten \calK$.
Thus we define
\eq{}
\tilPhi_{\calK}(A)=\tilp_{2!}\calF
\end{equation}
where $\tilp_2:~\calB\x X\to X$ is the projection to the
second multiple. It is again easy to see that the assignment 
$\calK\to \tilPhi_{\calK}$ reduces to a well-defined
morphism $\calK\to \tilphi_{\calK}$ from $K_{\calL^{-1}\boxtimes\calL,(w,w)}$ to 
$\Hom(K_{\calL,w},K_{\calL,w})$.

\prop{tilinversemap}
Let $\calK\in \calD_{\calL^{-1}\boxtimes \calL, (w,w)}(X\x X)$.
Then
\begin{enumerate}
\item $\tilphi_{\calK}(K^{\text{null}}_{\calL,w})\subset 
K^{\text{null}}_{\calL,w}$
\item Suppose that the image of $\calK$ in 
$K_{\calL^{-1}\boxtimes \calL,(w,w)}$ lies in 
$K^{\null}_{\calL^{-1}\boxtimes \calL,(w,w)}$. Then 
\eq{}
\tilphi_{\calK}(K_{\calL,w})\subset 
K^{\text{null}}_{\calL,w}
\end{equation}
It follows from the above that the correspondence $\calK\to \tilphi_{\calK}$
defines a map $\tilphi:~V_{\calL^{-1}\boxtimes\calL}\to \End (V_{\calL}$.
\item
$\tilphi\circ\psi=\id$
\item Let $\tilm(\calK)$ denote the modified matrix coefficient sheaf, 
corresponding
to $\calK$ (cf. \refss{modmat}) with the Weil structure, defined as in 
\refss{weil2}.
Then the $\qlb$-valued function $\tr(\tilm(\calK))$ is equal to the matrix 
coefficient
of the image of $\calK$ in $V_{\calL^{-1}\boxtimes \calL}=V_{\calL^{-1}}
\ten V_{\calL}$. 
\end{enumerate}
\eprop
The proof is analogous to the proof of Proposition \refp{inversemap}.
\ssec{}{Proof of Theorem \reft{main}$(2)$}
In order to show that the character of
$V_{\calL}$ is equal to $\tr (\chl)$, we must show that the image of 
$\al$ in $\Hom (V_{\calL},V_{\calL})$ is equal to the identity element. 
This follows easily from the following result.

\lem{} The functor $\tilPhi_{\al}$, defined by the sheaf $\al$, 
is canonically
isomorphic to identity functor. This isomorphism of functors 
commutes with the functor $\Phi_w\circ \fr^*$. 
\elem

The lemma is straightforward. On the other hand, together with
\refp{tilinversemap} it easily implies that
the character of $V_{\calL}$ is equal to $\tr(\chl)$, since it implies
that $\al$ represents the identity element in $V_{\calL^{-1}\boxtimes \calL}=
\Hom(V_{\calL},V_{\calL})$.

\sec{trace}{Trace formula}

Throughout this section
$X$ is a scheme of finite type over $k$,
where $k$ is an algebraically closed field.
For such a scheme we denote by $\D(X)$ the
bounded derived category of constructible $\qlb$-sheaves
(in \'etale topology) (see e.g. \cite{De}).
In this section we suggest a version of
Lefschetz-Verdier
trace formula in which the role of a geometric
correspondence is played by
an object of $\D(X\times X)$.
 
\subsection{Kernels and functors}
 
Recall that Verdier duality functor is
an equivalence of triangulated categories
$D:\D(X)^{op}\rightarrow\D(X):  A\mapsto
R\und{\Hom}(A, D_X)$ (where $D_X=p^!\qlb$,
$p$ is the projection to $\Spec(k)$),
such that $D^2\simeq\Id$,
 
For any $A,B\in\D(X)$ there is a natural isomorphism
\eq{RHom}
D(R\und{\Hom}(A,B))\simeq A\otimes DB.
\end{equation}
 
We consider functors $\Phi_K=\Phi_{K,!}:\D(X)\rightarrow\D(X)$
associated with kernels $K\in\D(X\times X)$:
\eq{phi}
\Phi_K(A)=p_{2!}(p_1^*A\otimes K).
\end{equation}
 
Another kind of functors from $\D(X)$ to itself is obtained
when composing $\Phi_K$ with Verdier duality. Namely, for
$K\in\D(X\times X)$ there is a functor
\eq{psi}
\Psi_K(A)=p_{2*}(R\und{\Hom}(K,p_1^!A))
\end{equation}
and the canonical isomorphism of functors
\begin{equation}
\Psi_K\simeq D\circ\Phi_K\circ D,
\end{equation}
derived from \refe{RHom} and standard isomorphisms
$D\circ p_{2!}\simeq p_{2*}\circ D$,
$D\circ p_1^*\simeq p_1^!\circ D$.
 
Let us introduce also the relative duality functor
$D_{p_1}:\D(X\times X)\rightarrow\D(X\times X)$ by setting
$D_{p_1}(K)=R\und{\Hom}(K,p_1^!\ov{\Q}_{l,X})$.
Note that we have a canonical morphism of functors
$\Id\rightarrow D_{p_1}^2$ which, however, is not an isomorphism
in general. Now we claim that
there is a canonical morphism of functors
\eq{fun_mor}
\Phi_K\rightarrow\Psi_{D_{p_1}(K)}
\end{equation}
 which is constructed as the
following composition
\begin{align*}
\Phi_K(A)=p_{2!}&(K\otimes p_1^*A)
\rightarrow \\
p_{2!}&(D_{p_1}^2(K)\otimes p_1^*A)
\wt{\rightarrow} p_{2!}(R\und{\Hom}(D_{p_1}(K),p_1^!\ov{\Q}_{l,X})\otimes
 p_1^*A)
\rightarrow \\
&\rightarrow p_{2!}(R\und{\Hom}(D_{p_1}(K),p_1^!A))=\Psi_{D_{p_1}(K)}(A)
\end{align*}
here we used the canonical morphism
$p_1^!\ov{\Q}_{l,X}\otimes p_1^*A\rightarrow p_1^!A$ which corresponds to
the natural morphism
$p_1^*A\rightarrow R\und{\Hom}(p_1^!\ov{\Q}_{l,X},p_1^!A)$.

\subsection{Trace functions}
 
With every $K\in\D(X\times X)$ we associate a finite-dimensional
vector space $V_K=\Hom(K,\De_*D_X)$ where $\De:X\rightarrow X\times X$
is the diagonal embedding. Then there is a canonical
isomorphism
\eq{homis}
V_K\simeq H^0_c(X,\De^*K)^*,
\end{equation}
which is constructed as follows:
$$
V_K=
\Hom(K,\De_*D_X)\simeq\Hom(\De^*K,D_X)\simeq H^0(X,D\De^*K)\simeq
H^0_c(X,\De^*K)^*
$$
where the last isomorphism is given by Verdier duality.
 
Let $\Phi_K:\D(X)\rightarrow\D(X)$ be a functor defined by \refe{phi},
$A$ be an object of $\D(X)$ equipped with a morphism
$\a:\Phi_K A\rightarrow A$. Then we can associate with
$\a$ an element $t_{K,A}(\a)\in V_K$ ("trace function") as follows.
By adjunction $\a$ corresponds to a morphism
$K\otimes p_1^*A\rightarrow p_2^!A$, or equivalently to a morphism
$\a':K\rightarrow R\und{\Hom}(p_1^*A, p_2^!A)$. Recall that
for every $\calF,\calG\in\D(X)$ there is
a canonical isomorphism
\eq{exterior}
R\und{\Hom}(p_1^*\calF, p_2^!\calG)\simeq D\calF\boxtimes \calG
\end{equation}
established in \cite{SGA5} (3.2). This isomorphism can be obtained
from the isomorphism $p_2^!\calG\simeq D_X\boxtimes \calG$ (which is the
particular case of \refe{exterior} for $\calF=\ov{\Q}_{l,X}$) using
the commutation of $R\und{\Hom}$ with the exterior tensor product.
Now we apply \refe{exterior} for $\calF=\calG=A$ and take the composition
of
$\a'$ with the natural morphism
$DA\boxtimes A\rightarrow\De_*(DA\otimes A)\rightarrow\De_*D_X$ to get an
element $t_{K,A}(\a)\in\Hom(K,\De_*D_X)= V_K$.
 
By duality, if we have an object $B\in\D(X)$ equipped with a
morphism $\b:B\rightarrow\Psi_K(B)$ then we obtain a
morphism $\b':\Phi_K(DB)\simeq D\Psi_K B\rightarrow DB$, hence, the above
construction
gives an element $t_{K,DB}(\b')\in V_K$.
Now if we are given an object $A\in\D(X)$ equipped with
a morphism $\a:A\rightarrow\Phi_K(A)$ we can take the composition of
$\a$ with the morphism \refe{fun_mor}
$\Phi_K(A)\rightarrow\Psi_{D_{p_1}(K)}(A)$
to get a morphism $A\rightarrow\Psi_{D_{p_1}}(A)$. Applying the above
remark we get an element of $V_{D_{p_1}(K)}$ which we denote by
$s_{K,A}(\a)$.
 
Note that we have a canonical isomorphism $p_1^!\ov{\Q}_{l,X}\simeq
p_2^*D_X$.
Hence, we get a canonical morphism
\begin{eqnarray}
d_K:\De^*D_{p_1}K\wt{\rightarrow}\De^*R\und{\Hom}(K,p_2^*D_X)\rightarrow
R\und{\Hom}(\De^*K,\De^*p_2^*D_X) \nonumber\\
\simeq R\und{\Hom}(\De^*K,D_X)=D\De^*K.
\end{eqnarray}
 
\defe{} An object $K\in\D(X\times X)$
is  called  {\it admissible}
if the composition of natural maps
$$
d_{K,*}:H^0_c(X,\De^*D_{p_1}K)\rightarrow H^0_c(X,D\De^*K)
\rightarrow H^0(X,D\De^*K)
$$
where the first arrow is induced by $d_K$, is an isomorphism.
An object $K$ is called {\it strictly admissible} if $d_K$ is an isomorphism
and the natural arrow $H^0_c(X,D\De^*K)\rightarrow H^0(X,D\De^*K)$ is an
isomorphism.
\edefe
 
\begin{exs} 1. If $f:B\rightarrow X\times X$ is a correspondence
with isolated fixed points, such that
$p_1\circ f$ is \'etale,
then $K=f_*(L)$ is strictly admissible for any local system $L$ on $B$.
 
\noindent
2. If $X$ is proper and $d_K$ is an isomorphism
then $K$ is strictly admissible. For example, $d_K$ is an isomorphism
for
$K=L\otimes p_2^*K'$ where $L$ is a local system on $X^2$, $K'\in\D(X)$.
\end{exs}
 
If $K$ is admissible then $d_{K,*}^{-1}$ gives an isomorphism
$$V_K\simeq H^0(X,D\De^*K)\wt{\rightarrow}
H^0_c(X,\De^*D_{p_1}K)
$$
Hence, we obtain a natural pairing
$$
\langle \cdot,\cdot \rangle: V_K\otimes V_{D_{p_1}K}\wt{\rightarrow}
H^0_c(X,\De^*D_{p_1}K)\otimes\Hom(\De^*D_{p_1}K, D_X)\rightarrow H^0_c(X,D_X)
$$
which induces a perfect pairing
$$
\Tr_X\langle \cdot,\cdot \rangle: V_K\otimes V_{D_{p_1}K} \rightarrow \qlb
$$
by composition with the trace map $\Tr_X:H^0_c(X,D_X)\rightarrow \qlb$.
 
{\bf Example.}\quad Let $f:B\rightarrow X\times X$ be a correspondence with
isolated fixed points such that $p_1\circ f$ is \'etale,
$K=f_*(\ov{\Q}_{l,B})$.
Then we have $D_{p_1}(K)\simeq K$.
Thus, we have a (symmetric) non-degenerate form $\chi$ on
$V_K\simeq H^0(B\cap\De,\qlb)^*$, hence on
$H^0(B\cap\De,\qlb)$,
where $B\cap\De:=B\times_{X\times X}\De$. The canonical decomposition
$H^0(B\cap\De,\qlb)\simeq
\oplus_{x\in B\cap\De}\qlb e_x$
is orthogonal with respect to this form, and for any $x\in B\cap\De$
we have $\chi(e_x,e_x)=m_x$,
where $m_x$ is the multiplicity of $x$ in the intersection-product
$B\cdot \De$.
Assume in addition
that $p_2\circ  f$  is  \'etale.  Then  the  trace  function
$t_{K,A}(\a)$
defined above for an object $A\in\D(X)$ and a morphism $\a:\Phi(A)\rightarrow
A$ is given by $t_{K,A}(\a)(e_x)=\Tr(\a_x,A_{\bar{x}})$ where
$\bar{x}=p_1f(x)=p_2f(x)\in X$, the endomorphism
$\a_x$ of $A_{\bar{x}}$ is the restriction to $x\in B$ of the morphism
$(p_1f)^*A\rightarrow  (p_2f)^!A\simeq  (p_2f)^*A$  corresponding to
$\a$.

\ssec{formula}{Formula}
 
Let $A,B\in\D(X)$ be a pair of objects equipped with
morphisms $\a:A\rightarrow\Phi_K(A)$ and $\b:\Phi_K(B)\rightarrow B$.
Then we have an induced endomorphism
of graded vector spaces:
\eq{phiab}
\phi(\a,\b):R\Hom(A,B)\rightarrow R\Hom(\Phi_K(A),\Phi_K(B))\rightarrow
R\Hom(A,B).
\end{equation}
Let $\phi_n(\a,\b):\Hom^n(A,B)\rightarrow\Hom^n(A,B)$ be the induced
morphisms. Then following the usual convention we denote
\begin{equation}
\Tr(\phi(\a,\b))=
\sum_i (-1)^i\Tr(\phi_i(\a,\b),\Hom^i(A,B)).
\end{equation}
 
On the other hand, we have defined trace functions
$s_{K,A}(\a)\in V_{D_{p_1}K}$ and $t_{K,B}(\b)\in V_K$.
Thus, if $K$ is admissible we can define the scalar
product $\Tr_X\lan s_{K,A}(\a),t_{K,B}(\b)\ran\in \qlb$.

\vspace{3mm}
\noindent
{\bf Conjecture}. {\it
Assume that $X$ is proper and $K$ is admissible. Then \it}
\eq{trace-formula}
\Tr_X\lan s_{K,A}(\a),t_{K,B}(\b)\ran = \Tr(\phi(\a,\b)).
\end{equation}

Below we will prove this under some additional technical assumptions.
One of these assumptions which we were unable to check has to do 
with some functoriality of the construction of the trace map.
Namely, for every admissible kernel $L$ and an object $E$ 
in $\D(X)$ we have defined the trace map 
$$s_{L,E}:\Hom(p_2^*E, p_1^*E\otimes L)\simeq\Hom(E,\Phi_L(C))\rightarrow 
V_{D_{p_1(L)}}\simeq H^0(X,\De^*L)$$ 
(the first isomorphism here is due to the fact that $X$ is proper). 
Now let $f:X'\rightarrow X$ be a proper morphism such that $(f\times f)^*L$ 
is admissible. Then we'll say that the triple $(f,L,E)$ {\it behaves
functorially} if the following diagram commutes: 
\eq{com1}
\setlength{\unitlength}{0.30mm} 
\begin{array}{ccc} 
\Hom(p_2^*E,p_1^*E\otimes L) & \lrar{s_{L,E}} & H^0(X,\De^*L) \\ 
\ldar{} & & \ldar{} \\ 
\Hom(p_2^*f^*E,p_1^*f^*E\otimes(f\times f)^* L) & 
\lrar{s_{(f\times f)^*L,f^*E}} & H^0(X',f^*\De^*L). 
\end{array} 
\end{equation}

\th{trace-theorem}   Assume  that  $X$  is  proper,  $K$  is
admissible, the triple 
$$(\De:X\rightarrow X\times X, L=p_{13}^*K\otimes p_{24}^*D_{p_1}K, 
E=A\boxtimes\ov{B})$$
behaves functorially, and the morphism $\b$ factorizes as follows
$$\b:\Phi_K(B)\rightarrow\Psi_{D_{p_1}(K)}(B)\stackrel{\wt{\b}}{\rightarrow} 
B,$$
where the first arrow is given by the canonical morphism
\refe{fun_mor}.
Then the formula \refe{trace-formula} holds.
\eth

\smallskip
                       
\Pf .
The idea of the proof is to reduce \refe{trace-formula}
to the usual Lefschetz-Verdier formula.
By adjointness and duality the pair of
morphisms $\a$ and $\b$ corresponds
to a pair of morphisms $\a':p_2^*A\rightarrow p_1^*A\otimes K$ and
$\b':p_2^*\ov{B}\rightarrow R\und{\Hom}(K,p_1^!(\ov{B}))$ where
$\ov{B}=D(B)$.
Taking the tensor product of $\a'$ and $\b'$ we obtain the morphism
$$
p_2^*(A\otimes\ov{B})\rightarrow
p_1^*A\otimes K\otimes R\und{\Hom}(K,p_1^!(\ov{B}))\rightarrow
p_1^*A\otimes p_1^!\ov{B}\rightarrow p_1^!(A\otimes\ov{B}),
$$
where the last arrow is induced by
the canonical isomorphism $p_1^!E\simeq p_1^*E\otimes
p_2^*D_X$ for any $E\in \D(X)$.
Let $C=A\otimes\ov{B}$, $\ga:p_2^*C\rightarrow p_1^!C$
be the morphism defined above.
Applying Lefschetz-Verdier formula (Thm.~3.3 of \cite{SGA5})
to $\ga$ and the diagonal correspondence
$p_1^*C\rightarrow\De_*C\rightarrow p_2^!C$
we obtain the equality
\eq{Verdier-formula}
\Tr(\ga_*)=\Tr_X(\ga_{\De})
\end{equation}
where the LHS is the trace of the induced map
$\ga_*: R\Ga(X,C)\rightarrow R\Ga(X,C)$. The RHS is obtained by
applying the the trace map
$\Tr_X:H^0_c(X,D_X)\stackrel{\sim}{\rightarrow}\qlb$ to the morphism
$\ga_{\De}:\ov{\Q}_{l,X}\rightarrow D_X$ obtained
by adjointness from the following morphism induced by $\ga$:
$$
\ov{\Q}_{l,X^2}\rightarrow R\und{\Hom}(p_2^*C,p_1^!C)
\simeq C\boxtimes DC\rightarrow
\De_*(C\otimes DC)\rightarrow\De_*D_X
$$
It follows from Lemmas \refl{glob} and \refl{loc} below
that $\Tr(\phi(\a,\b))=\Tr(\ga_*)$ and
$\Tr_X(\ga_{\De})=\Tr_X\lan s_{K,A}(\a),t_{K,B}(\b)\ran$.
Hence, \refe{trace-formula} follows from \refe{Verdier-formula}.
\ed
 
\lem{glob} Let $X$ be proper. Then
under the  isomorphism
$$
R\Hom(A,B)^*\simeq R\Ga(X,D(R\und{\Hom}(A,B)))\simeq R\Ga(X,A\otimes \ov{B})
$$
one has $\phi(\a,\b)^*=\ga_*$ where $\phi(\a,\b)^*$ is the dual
operator to \refe{phiab}.
\elem
 
\Pf  By definition $\phi=\phi(\a,\b)$ is the composition of
the natural maps
$$
f:R\Hom(A,B)\rightarrow R\Hom(p_1^*A\otimes K, p_1^*B\otimes K)
$$
and
$$
g:R\Hom(p_1^*A\otimes K, p_1^*B\otimes K)
\rightarrow R\Hom(\Phi_K(A),\Phi_K(B))
$$
and the map
$$
h=h(\a,\b):R\Hom(\Phi_K(A),\Phi_K(B))\rightarrow R\Hom(A,B)
$$
induced by $\a$ and $\b$. We have the following natural isomorphisms:
\begin{align*}
&R\Hom(A,B)^*\simeq R\Ga(X,A\otimes\ov{B})\\
&R\Hom(p_1^*A\otimes K,p_1^*B\otimes K)^*\simeq R\Ga(X\times X,
p_1^*A\otimes K\otimes D(p_1^*B\otimes K))\simeq\\
&\simeq R\Ga(X\times X,
p_1^*A\otimes K\otimes R\und{\Hom}(K,p_1^!\ov{B})),\\
&R\Hom(\Phi_K(A),\Phi_K(B))^*\simeq R\Ga(X,\Phi_K(A)\otimes D(\Phi_K(B)))
\simeq R\Ga(X,\Phi_K(A)\otimes \Psi_K(\ov{B})).
\end{align*}
Under these identifications $f^*$ is induced by the canonical morphism
\eq{f*}
p_{1!}(p_1^*A\otimes K\otimes R\und{\Hom}(K,p_1^!\ov{B}))\rightarrow
p_{1!}(p_1^*A\otimes p_1^!\ov{B})\rightarrow p_{1!}p_1^!(A\otimes\ov{B})
\rightarrow A\otimes\ov{B},
\end{equation}
$g^*$ is induced by the morphism
\begin{align*}
\Phi_K(A)\otimes\Psi_K(\ov{B})&\rightarrow p_{2^*}(p_1^*A\otimes K)\otimes
p_{2*}(R\und{\Hom}(K,p_1^!\ov{B}))\rightarrow \\
&p_{2*}(p_1^*A\otimes K\otimes R\und{\Hom}(K,p_1^!\ov{B}))
\end{align*}
and $h$ is induced by the morphism
$$
\a\otimes D(\b):A\otimes\ov{B}\rightarrow\Phi_K(A)\otimes\Psi_K(\ov{B})
$$
It follows that $g^*h^*$ is induced by the morphism
\eq{g*h*}
A\otimes\ov{B}\rightarrow p_{2*}p_2^*(A\otimes\ov{B})
\overset{p_{2*}(\a'\otimes\b')}\to
p_{2*}(p_1^*A\otimes K\otimes R\und{\Hom}(K,p_1^!\ov{B}))
\end{equation}
and comparing \refe{f*} and \refe{g*h*}
with the definition of $\ga:p_2^*C\rightarrow p_1^!C$
(where $C=A\otimes\ov{B}$) we conclude that
$\phi^*=f^*g^*h^*$ is equal to the following composition
\eq{phi*}
\begin{align*}
R\Ga(X,C)\rightarrow &R\Ga(X\times X,p_2^*C)\rightarrow\\ 
&R\Ga(X, p_{1!}p_2^*C)
\overset{R\Ga(p_{1!}
(\ga))}{\to}
R\Ga(X,p_{1!}p_1^!C)\rightarrow R\Ga(X,C).
\end{align*}
\end{equation}
But this is the definition of $\ga_*$.
\ed
 
\lem{loc} Under the assumptions of the theorem one has
\eq{}
\ga_{\De}=\lan s_{K,A}(\a),t_{K,B}(\b)\ran
\end{equation}
\elem
 
\prf  Let us apply the functoriality assumption. 
We  have  a   morphism   $\a':p_2^*A\to   p_1^*A\otimes   K$
corresponding to $\a$
and a morphism $\wt{\b'}:p_2^*\ov{B}\to p_1^*\ov{B}\otimes
D_{p_1}K$  corresponding to $\wt{\b}$. Their external
product is the morphism
$$
\a'\boxtimes\wt{\b'}:p_3^*A\otimes p_4^*\ov{B}\to
p_1^*A\otimes p_2^*\ov{B}\otimes p_{13}^*K\otimes
p_{24}^*D_{p_1}K.
$$
Its trace $s_{L,E}(\a'\boxtimes\wt{\b'})\in H^0(X\times X,
\De^*K\boxtimes\De^*D_{p_1}K)$ is equal to
the external product of $s_{K,A}(\a)\in H^0(X,\De^*K)$
and $t_{K,B}(\b)\in H^0(X,\De^*D_{p_1}K)$.
Hence, the commutative diagram \refe{com1} tells us in this
case that the element
$$
s_{K,A}(\a)\cup t_{K,B}(\b)=
\De^*s_{L,E}(\a'\boxtimes\wt{\b'})\in H^0(X,\De^*(K\otimes
D_{p_1}K))
$$
is obtained as the trace of the morphism
$$
\a'\otimes\wt{\b'}:
p_2^*(A\otimes\ov{B})\to p_1^*(A\otimes\ov{B})\otimes
(K\otimes D_{p_1}K).
$$
The scalar product $\lan s_{K,A}(\a), t_{K,B}(\b)\ran$
is the image of $s_{K,A}(\a)\cup t_{K,B}(\b)$
under the map $H^0(X,\De^*(K\otimes D_{p_1}K))\to H^0(X,D_X)$
induced by the natural morphism
$K\otimes D_{p_1}K\to p_2^*D_X$.
On the other hand, composing $\a'\otimes\wt{\b'}$ with the
latter morphism we obtain the morphism
$$
\ga:p_2^*(A\otimes\ov{B})\to p_1^*(A\otimes\ov{B})\otimes
p_2^*D_X
$$
introduced above. Hence, its trace $\ga_{\De}$ is equal to
$\lan s_{K,A}(\a), t_{K,B}(\b)\ran$ as required.
\epr

\rem{}
For non-proper $X$ the formula \refe{trace-formula} holds also for
the Frobenius correspondence.
One can ask by analogy with Deligne conjecture, whether
\refe{trace-formula} holds for a composition of a given functor
with a sufficiently high power of Frobenius correspondence.
For example, one can show that it holds for the composition of
the symplectic Fourier-Deligne transform
with Frobenius correspondence.
\erem

\sec{delus}{Deligne-Lusztig versus Kazhdan-Laumon representations}

In this section we will assume that formula \refe{trace-formula} holds
(in fact, we will make a slightly different assumption, whose proof, however,
should be the same as that formula  \refe{trace-formula}).
Modulo this assumption we will explain how to connect geometrically
the Kazhdan-Laumon representation with Deligne-Lusztig representation
(\cite{dl}).

\ssec{}{Deligne-Lusztig representations}In this subsection we review 
the definition and basic properties of the Deligne-Lusztig representations
(cf. \cite{dl}). Our notations, however, will be different from \cite{dl}.
\sssec{}{The varieties $X_w$}Recall that in \refsss{Z-w} we have
defined for every $w\in W$ certain smooth closed subvariety $Z_w$ 
inside $X\x X$ of dimension $\dim(X)+l(w)$. Moreover, according to 
Proposition \refsss{Z-w}(4), we know that  
that $Z_w$ intersects $\Gam_{\fr}$
transversally, where $\Gam_{\fr}$ denotes the graph of the Frobenius
morphism on $X$. 

Set now $X_w=Z_w\cap \Gam_{\fr}$. Then it follows from
the above that $X_w$ is a reduced smooth closed subvariety in $X\x X$,
which can be also regarded as a closed subvariety of $X$ using the
projection on the first factor. One can easily see that $\dim(X_w)=l(w)$. 
$X_w$ admits a natural action of the
group $\gfq\x T_w(\fq)$, which is inherited from its action on $X$.
\sssec{delus-theorem}{The representations $R_{\theta,w}$}
Consider now $H^*_c(X_w,\qlb)$ -- the $\ell$-adic cohomology of
$X_w$ with compact supports. Since the finite group $\gfq\x T_w(\fq)$
acts on $X_w$ it acts also on $H^*_c(X_w,\qlb)$.

Thus we may
consider $\sum (-1)^iH^i_c(X_w,\qlb)$ as a virtual representation
of $\gfq$ and decompose it with respect to characters of $T_w(\fq)$.
For any character $\theta:\ T_w(\fq)\to \qlb^*$ we will denote
by $R_{{\theta},w}$ the corresponding virtual representation
of $G(\fq)$. More precisely, for any $i\geq 0$ let 
\eq{}
H^i_c(X_w,\qlb)_{\theta}=\{ \xi\in H^i_c(X_w,\qlb)|\ t(\xi)=\theta(t)\xi
\quad\text{for any $\xi\in T_w(\fq)$}\}
\end{equation}
and we define 
\eq{}
R_{\theta,w}=\sum\limits_i (-1)^iH^i_c(X_w,\qlb)_{\theta^{-1}}
\end{equation}
The next statement is due to P.~Deligne and G.~Lusztig (cf. \cite{dl})
for the case when $q$ is large enough, and it is due to 
B.~Haastert (cf. \cite{ha}) in the general case. 
\th{delus-theorem} 
Suppose that $\theta$ is quasi-regular (cf. \refss{characters}). Then

1) the natural map of forgetting the supports from 
$H^*_c(X_w,\qlb)_{\theta}$ to $H^*(X_w,\qlb)_{\theta}$
is an isomorphism and $H^i_c(X_w,\qlb)=0$ for $i\neq l(w)$.

2) One has 
\eq{hom-delus}
\dim\Hom_{\gfq}(H^{l(w)}_c(X_w,\qlb),H^{l(w)}_c(X_w,\qlb))=\#\wtheta^{\frw}
\end{equation}
In particular, $(-1)^{l(w)}R_{\theta,w}$ is an irreducible representation
of $\gfq$ if $\theta$ is regular.
\eth
\ssec{}{The case of quasi-regular $\calL$}In this subsection we suppose that
$\calL$ is non-singular. 
We assume that an analogue of the
formula \refe{trace-formula} holds when $X$ is the basic affine space of
$G$, $K=(\id\times\fr)*\QQ_{\ell,Z_w}[l(w)]$, 
$(A,\a)\in\Perv_{\calL,w}(X)$ and $(B,\b)\in\Perv_{\calL,w}(X)$.
More precisely, in this case 
we can define analogues of LHS and RHS of \refe{trace-formula}
as follows.
First we notice that $\Phi_K(A)\simeq\Phi_w(\fr^*A)$
and that the morphism
$d_K:\De^*D_{p_1}K\rightarrow D\De^*K$ is an isomorphism in our case.
Thus, we can consider the trace function
$s_{K,A}(\a)\in V_{D_{p_1}K}\simeq H^0(X,\De^*K)$.
The diagonal action of $T_w(\fq)=T(\ov{\fr})^{\fr_w}$ on $X\times X$
lifts to an action on $K$, such that the induced action of
$T_w(\fq)$ on $H^0(X,\De^*K)\simeq H^{l(w)}(X_w,\qlb)$
is the natural one. 
It is easy to see that 
$s_{K,A}(\a)$ lies in $\theta^{-1}$-component
of $H^0(X,\De^*K)$ where $\theta=\tr\calL$ is
the corresponding character of $T_w(\fq)$.
Similarly, we can consider the trace $t_{K,B}(\b^{-1})$ of 
$\b^{-1}:\Phi_K(B)\rightarrow B$ which lies in
$H^0(X,\De^*D_{p_1}K)_{\theta}$.
Now we have $\De^*K\simeq \QQ_{\ell,X_w}[l(w)]$, hence,
by theorem \reft{delus-theorem} the natural map
$H^0_c(X,\De^*K)_{\theta}\rightarrow H^0(X,\De^*K)_{\theta}$
is an isomorphism. In particular, the intersection pairing
gives rise to the pairing
$$H^0(X,\De^*K)_{\theta^{-1}}\times H^0(X,\De^*D_{p_1}K)_{\theta}
\rightarrow\qlb$$
so we can define $\lan s_{K,A}(\a), t_{K,B}(\b^{-1})\ran\in\qlb$.

To define the analogue of RHS of  \refe{trace-formula}
we remark that in the proper case we used the traces of
endomorphisms induced by $\a$ and $\b$ on hypercohomologies of
$A\otimes D(B)$. In our case $A\otimes D(B)$ is equivariant
with respect to the action of $T$ on $X$, hence it descends
to a sheaf $\ov{A\otimes D(B)}$ on the flag variety $X/T$.
Furthermore, since $\a$ and $\b$ were compatible with the action
of $T$ we get the induced endomorphisms $\ov{\phi}(\a,\b)$
on hypercohomologies
of $\ov{A\otimes D(B)}$ so we can take the alternated sum of their
traces. The obtained number differs from the pairing
$\<(A,\a),(D(B),D(\b)^{-1})\>$ defined in \refsss{w_arb} 
by a constant non-zero multiple.

Next we observe that the morphism of functors
$\Phi_K(B)\rightarrow\Psi_{D_{p_1}K}(B)$ is an isomorphism, 
so the last condition of theorem \reft{trace-theorem}
is satisfied. Thus, assuming the functoriality property
in the formulation of this theorem we can slightly
modify the argument to prove the equality  
\eq{trace-formula2}
\Tr_X\lan s_{K,A}(\a),t_{K,B}(\b^{-1})\ran = \Tr(\ov{\phi}(\a,\b)).
\end{equation}

The following result is
obtained as a corollary of the formula \refe{trace-formula2}. 

\th{} Let $\calL$ be a quasi-regular one-dimensional local system on $T$ and
let $w\in W$ be an element of the Weyl group, such that $\frw^*(\calL)$ is
isomorphic to $\calL$. Let $\theta$ be the corresponding character of 
$T_w(\fq)$. Then one has canonical isomorphism of $\gfq$-representations
\eq{isomorphism}
V_{\calL}\simeq H^{l(w)}_c(X_w)_{\theta}^*\simeq R_{\theta,w}
\end{equation}
\eth

\prf The proof of the second isomorphism follows merely from 
Theorem \reft{delus-theorem}. Therefore, it is enough to construct
the first isomorphism. This is done in the following way.

We have to define a map $\kappa:\ V_{\calL}\to H^{l(w)}_c(X_w)_{\theta}^*$. 
We will do
it as follows. First, we will define $\okap:K_{\calL,w}\to 
H^{l(w)}_c(X_w)_{\theta}^*$
(which will be a map of $\gfq$-modules) and then show that it is surjective
and that its kernel is equal to $K^{\text{null}}_{\calL,w}$. 

So, we define
\eq{}
\okap((A,\alp))=t_{A,\alp}
\end{equation}
for any $(A,\alp)\in \Perv_{\calL,w}(X)$. It is easy to see that $\okap$
is a well defined map from $K_{\calL,w}$ to $H^{l(w)}_c(X_w)_{\theta}^*$, which
commutes with $\gfq$-action. 
We claim, first of all, that $K^{\text{null}}_{\calL,w}$
contains the kernel of $\okap$. Indeed, suppose that we
have some $(A,\alp)\in \Perv_{\calL,w}(X)$ whose image in $K_{\calL,w}$
does not lie in $K^{\text{null}}_{\calL,w}$. Then there exists
some $(B,\beta)\in \Perv_{\calL^{-1},w}$, such that
$\<(A,\alp),(B,\beta)\>\neq 0$. It follows from 
formula \refe{trace-formula2} (which we assumed to hold in this case),
that $\<t_{A,\alp},s_{D(B),D(\beta)^{-1}}\>$ is also non-zero.
Hence, $\<t_{A,\alp},s_{D(B),D(\beta)^{-1}}\>\neq 0$ and 
therefore $t_{A,\alp}\neq 0$,
which is what we had to prove.

It follows now from the fact that 
$\Ker\okap\subset K^{\text{null}}_{\calL,w}$ that the map $\okap$ identifies 
$V_{\calL}$ with a subquotient of $H^{l(w)}_c(X,\qlb)_{\theta}^*$. Therefore,
in order to show that $\okap$ in fact descends to an isomorphism between the
two, it is enough to show (since $\gfq$ is a finite group) that 
\eq{hom-equality}
\dim\Hom_{\gfq}(V_{\calL},V_{\calL})=
\dim\Hom_{\gfq}(H^{l(w)}_c(X,\qlb)_\theta^*,H^{l(w)}_c(X,\qlb)_\theta^*)
\end{equation}

(In fact, we already know, in fact, that the two representations are 
isomorphic, since the character of both of them is equal to $\tr(\chl)$.
However, we want to show that $\okap$ defines an isomorphism between
$V_{\calL}$ and $(-1)^{l(w)}R_{\theta,w}$ independently of the
above computation of their characters). 

We know now that the left hand
side of \refe{hom-equality} is equal to $\#\wl^{\frw}$ and the right hand
side is equal to $\#\wtheta^{\frw}$. On the other hand, it is easy
to see that $\wl^{\frw}=\wtheta^{\frw}$, which finishes the proof.
\epr



\begin{thebibliography}{BGG2}

\bibitem{bgg} A.~Beilinson, J.~Bernstein and P.~Deligne,
{\em Faisceaux pervers}, Ast\'erisque, {\bf 100} (1982) 

\bibitem{bl} J.~Bernstein and V.~Lunts,
{\em Equivariant sheaves and functors}, 
Lecture Notes in Math.

\bibitem{br} J.-L.~Brylinski, 
{\em Transformations canoniques, dualit\'e projective, the\'eorie de
Lefschetz, transformations de Fourier et sommes trigonom\'etriques},
Ast\'erisque {\bf 140} (1986), 3-134

\bibitem{De} P.~Deligne, {\it La conjecture de Weil II},
Publ. I.H.E.S.,  {\bf 52} (1980), 137-252 

\bibitem{De1} P.~Deligne, {\it Action du groupe des tresses sur une
cat\'egorie}, Invent.
   Math. {\bf 128} (1997), no. 1, 159--175

\bibitem{dl}
P.~Deligne and G.~Lusztig, 
{\em Representations of reductive groups over a finite field},
Ann. of Math., {\bf 103}(1976) 103-161

\bibitem{fu}
K.~Fujiwara,
{\em Rigid geometry, Lefschetz-Verdier trace formula and 
Deligne's conjecture}, 
Invent. Math. {\bf 127} (1997) 489--533

\bibitem{ha}
R.~Haastert, 
{\em Die Quasiaffinit\"at der Deligne-Lusztig Variet\"aten}, 
J. of Alg. {\bf 102} (1986), 186-193

\bibitem{KatzLaum}
N.~Katz and G.~Laumon,
{\em Transformation de Fourier et majoration de
sommes exponentielles}, 
Publ. IHES, {\bf 62} (1985), 361--418.

\bibitem{ka}
D.~Kazhdan,
{\em ``Forms'' of the principle series for $GL_n$},
in: {\em Functional analysis on the Eve of the 21-st century},
Progress in Math., Birkh\"auser, {\bf 131} (1995) 153-172
 
\bibitem{kl}
D.~Kazhdan and G.~Laumon,
{\em Gluing of perverse sheaves and discrete series representations},
Jour. of Geom. and Physics, {\bf 5} (1988), 63-120

\bibitem{laum}
G.~Laumon,
{\em Faisceaux characters (d'apr\`es Lusztig)},
S\'em. Bourbaki, Ast\'erisque, {\bf 178-179} (1989), 231-260

\bibitem{lu-book}
G.~Lusztig, 
{Characters of reductive groups over a finite field},
Ann. of Math. Stud., 107, Princeton Univ. Press,
Princeton, NJ
(1984)

\bibitem{lu-char}
G.~Lusztig
{\em Character sheaves I}, Adv. in Math. {\bf 56} (1985), 193-237

 
\noindent
G.~Lusztig
{\em Character sheaves II}, Adv. in Math. {\bf 57} (1985), 226-265

G.~Lusztig
{\em Character sheaves III}, Adv. in Math. {\bf 57} (1985), 266-315

G.~Lusztig
{\em Character sheaves IV}, Adv. in Math. {\bf 59} (1986), 1-63

G.~Lusztig
{\em Character sheaves V}, Adv. in Math. {\bf 61} (1986), 103-155

\bibitem{lu-86}
G.~Lusztig,
{\em On the character values of finite 
Chevalley groups at unipotent elements},
Journal of Algebra {\bf 104} (1986), 146-194

\bibitem{lu-obzor}
G.~Lusztig,
{\em Introduction to character sheaves},
Proceedings of Symposia in Pure Math. {\bf 47} (1987), 165-179

\bibitem{lu-90}
G.~Lusztig,
{\em Green functions and character sheaves},
Annals of Math., {\bf 131}, 1990, 355-408

\bibitem{mir-vil}
I.~Mirkovic and K.~Vilonen,
{\em Characteristic varieties of character sheaves},
Inventiones Math. {\bf 93} (1988), 405-418

\bibitem{po} A.~Polishchuk, 
{\em }, Preprint (1997).

\bibitem{SGA4} {\em Th\'eorie de topos et cohomologie des
sch\`emas}, Tome 3, SGA 4, Lecture Notes in Mathematics,
{\bf 305} (1973)

\bibitem{SGA5}
{\em Cohomologie $\ell$-adique et fonctions $L$}, SGA 5, 
Edit\'e par L.~Illusie,
Lecture Notes in Mathematics,{\bf 589}, 1977

\end{thebibliography}
\end{document}